\documentclass[11pt]{article}
\usepackage[utf8]{inputenc}
\usepackage{amsmath,amssymb,epsfig,bbm}
\usepackage{stmaryrd}
\usepackage{comment}
\usepackage{color}
\usepackage[T1]{fontenc}

\usepackage[textsize=small]{todonotes}
\usepackage{enumitem}
\usepackage{varwidth}
\setlist{nolistsep}
\usepackage{hyperref}

\newtheorem{defi}{Definition}
\newtheorem{prop}[defi]{Proposition}
\newtheorem{theo}[defi]{Theorem}
\newtheorem{conj}[defi]{Conjecture}
\newtheorem{lemm}[defi]{Lemma}
\newtheorem{coro}[defi]{Corollary}
\newtheorem{rema}[defi]{Remark}
\newtheorem{exem}[defi]{Example}
\newtheorem{exems}[defi]{Examples}

\newcommand{\bdefi}{\begin{defi}}
\newcommand{\edefi}{\end{defi}}
\newcommand{\bprop}{\begin{prop}}
\newcommand{\eprop}{\end{prop}}
\newcommand{\btheo}{\begin{theo}}
\newcommand{\etheo}{\end{theo}}
\newcommand{\blemm}{\begin{lemm}}
\newcommand{\brema}{\begin{rema}}
\newcommand{\erema}{\end{rema}}
\newcommand{\bexer}{\begin{exem}}
\newcommand{\eexer}{\end{exem}}
\newcommand{\bexems}{\begin{exems}}
\newcommand{\eexems}{\end{exems}}
\newcommand{\bconj}{\begin{conj}}
\newcommand{\econj}{\end{conj}}
\newcommand{\elemm}{\end{lemm}}
\newcommand{\bcoro}{\begin{coro}}
\newcommand{\ecoro}{\end{coro}}
\newcommand{\dem}{\noindent{\bf Proof. }}
\newcommand{\rem}{\noindent{\bf Remark. }}

\usepackage{mathrsfs}
\renewcommand\mathcal{\mathscr}

\renewcommand{\H}{{\cal H}}

\newcommand{\M}{{\cal M}}
\newcommand{\N}{{\cal N}}
\newcommand{\OOO}{{\cal O}}

\newcommand{\maths}[1]{{\mathbb #1}}  

\newcommand{\BB}{\maths{B}}
\newcommand{\CC}{\maths{C}}

\newcommand{\FF}{\maths{F}}
\newcommand{\HH}{\maths{H}}
\newcommand{\KK}{\maths{K}}

\newcommand{\NN}{\maths{N}}
\newcommand{\OO}{\maths{O}}
\newcommand{\PP}{\maths{P}}
\newcommand{\QQ}{\maths{Q}}
\newcommand{\RR}{\maths{R}}
\newcommand{\SSS}{\maths{S}}

\newcommand{\UU}{\maths{U}}

\newcommand{\ZZ}{\maths{Z}}

\newcommand{\mmm}{{\mathfrak m}}

\newcommand{\weakstar}{\overset{*}\rightharpoonup}
\newcommand{\ra}{\rightarrow}
\newcommand{\bs}{\backslash}

\newcommand{\ov}[1]{{\overline #1}} 
\newcommand{\wt}[1]{{\widetilde{#1}}}

\newcommand{\ga}{\gamma}
\newcommand{\Ga}{\Gamma}

\newcommand{\cqfd}{\hfill$\Box$}

\newcommand{\bigO}{\operatorname{O}}
\newcommand{\card}{{\operatorname{Card}}}

\newcommand{\covol}{\operatorname{Covol}}

\newcommand{\diam}{{\operatorname{diam}}}

\newcommand{\dvol}{\;d\operatorname{vol}}
\newcommand{\Gal}{\operatorname{Gal}}

\newcommand{\haarheis}{\operatorname{Haar}_{\operatorname{Heis}_3}}
\newcommand{\Heis}{\operatorname{Heis}}
\newcommand{\HS}{\mathcal{H\!S}}
\newcommand{\id}{\operatorname{id}}
\renewcommand{\Im}{{\operatorname{Im}}}
\newcommand{\Isom}{\operatorname{Isom}}
\newcommand{\Leb}{\operatorname{Leb}}

\renewcommand{\Re}{{\operatorname{Re}}}

\newcommand{\smallo}{\operatorname{o}}

\newcommand{\Vol}{\operatorname{Vol}}
\newcommand{\vol}{\operatorname{vol}}

\newcommand{\hnr}{{\HH}^n_\RR}
\newcommand{\hnc}{{\HH}^n_\CC}
\newcommand{\hdc}{{\HH}^2_\CC}

\newcommand{\PSL}{\operatorname{PSL}}
\newcommand{\SL}{\operatorname{SL}}

\newcommand{\SU}{\operatorname{SU}}
\newcommand{\PSU}{\operatorname{PSU}}

\newcommand{\tr}{\operatorname{\tt tr}}
\newcommand{\n}{\operatorname{\tt n}}

\newcommand\normalout{\partial^1_{+}}

\newcommand\normalpm{\partial^1_{\pm}}

\newcounter{fig}

\title{Counting and equidistribution in Heisenberg groups}
\author{Jouni Parkkonen \and Fr\'ed\'eric Paulin} 

\begin{document}
\bibliographystyle{../alphanum}
\maketitle
\begin{abstract}
  We strongly develop the relationship between complex hyperbolic
  geometry and arithmetic counting or equidistribution applications,
  that arises from the action of arithmetic groups on complex
  hyperbolic spaces, especially in dimension $2$.  We prove a Mertens'
  formula for the integer points over a quadratic imaginary number
  fields $K$ in the light cone of Hermitian forms, as well as an
  equidistribution theorem of the set of rational points over $K$ in
  Heisenberg groups.  We give a counting formula for the cubic points
  over $K$ in the complex projective plane whose Galois conjugates are
  orthogonal and isotropic for a given Hermitian form over $K$, and a
  counting and equidistribution result for arithmetic chains in the
  Heisenberg group when their Cygan diameter tends to $0$.
\end{abstract}

\section{Introduction}

The aim of this paper is to give original asymptotic counting and
equidistribution results with error terms of arithmetically defined
points or circles in nilmanifolds covered by the Heisenberg groups. We
refer for instance to \cite{Breuillard05, BjoFis09, GreTao12, Kim13,
  BenQui13a} or \cite[Chap.~10]{EinWar11} for other types of results.

One of our main results (see Section \ref{sect:chains}) is an
asymptotic counting result of arithmetic chains in hyperspherical
geometry. Let $q$ be the Hermitian form $-z_0\overline{z_2}
-z_2\overline{z_0}+|z_1|^2$ on $\CC^3$. Its isotropic locus in the
complex projective plane $\PP_2(\CC)$ is called by Poincaré the {\it
  hypersphere} \cite{Poincare07}. With $[1:0:0]$ removed, the
hypersphere identifies with the $3$-dimensional Heisenberg group
$\Heis_{3}$ (central extension of $\CC$ by $\RR$), and it carries a
natural distance $d_{\rm Cyg}$, called Cygan's (sometimes Koranyi's)
distance. Recall that a {\it chain}, as introduced by von Staudt and
developped in particular by E.~Cartan, see for instance
\cite{Cartan32} and \cite[\S 4.3]{Goldman99}, is a nontrivial
intersection with the hypersphere of a complex projective line. A
chain is either a fiber of the canonical morphism $\Heis_{3} \ra\CC$
or an ellipse whose projection by this map is a circle. Let $K$ be an
imaginary quadratic number field, and let $\OOO_K$ be its ring of
integers. We say that a chain $C_0$ is {\it arithmetic} (over $K$) if
the orbit of some point in $C_0$ under the stabiliser of $C_0$ in the
arithmetic lattice $\PSU_q(\OOO_K)$ is dense in $C_0$. The stabiliser
$\PSU_q(\OOO_K)_\infty$ of $[1:0:0]$ in $\PSU_q(\OOO_K)$ preserves the
diameters of the chains for $d_{\rm Cyg}$. The picture below shows an
orbit of arithmetic chains under the arithmetic lattice
$\PSU_q(\ZZ[i])$.

\begin{center}
\includegraphics{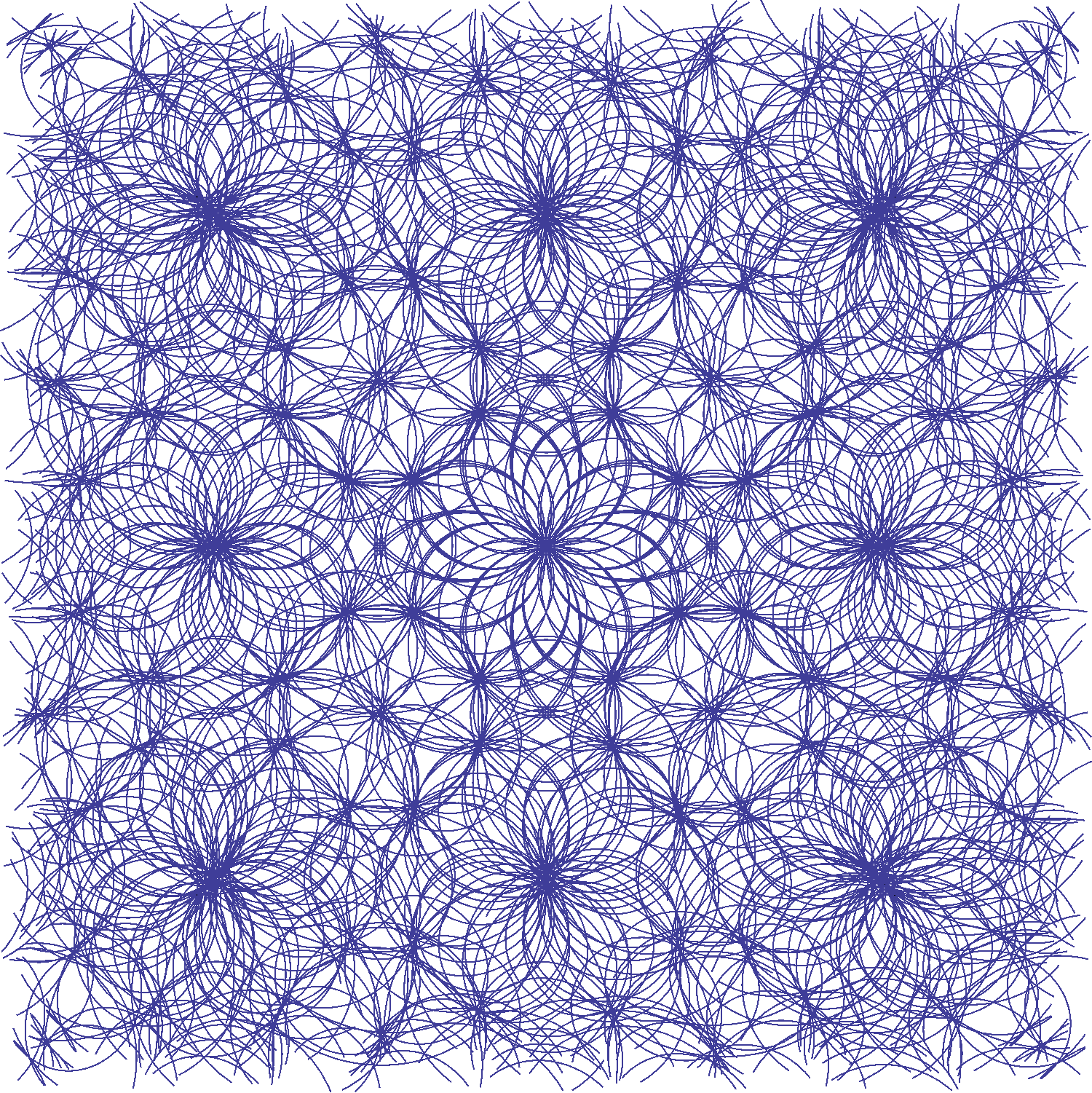}
\end{center}

\btheo\label{theo:countingchainintro} Let $C_0$ be an arithmetic chain
in the hypersphere.  There exists a constant $\kappa>0$ and an
explicit constant $c>0$ such that, as $s\ra+\infty$, the number of
chains modulo $\PSU_q(\OOO_K)_\infty$ in the $\PSU_q(\OOO_K)$-orbit of
$C_0$, with diameter at least $\epsilon$, is equal to
$c\;\epsilon^{-4}(1+\bigO(\epsilon^{\kappa}))$.  
\etheo

We will also prove that the centers of the finite arithmetic chains
equidistribute in the hypersphere.

An analogous method allows us in Section \ref{sect:cubicpoints} to
prove the following counting result of some arithmetic points in the
complex projective plane $\PP_2(\CC)$.  Let $z_0\in\PP_2(\CC)$ be a
cubic point over $K$, whose Galois conjugates $z'_0,z''_0$ over $K$
are isotropic and orthogonal to $z_0$ for the Hermitian form $q$. The
arithmetic lattice $\PSU_q(\OOO_K)$ acts with infinitely many orbits
on the set of such points.  The inverse of a slightly modified version $d''_{\rm Cyg}$ (see Section
\ref{sec:cxhyp} for precise definitions) of the Cygan distance between
the two isotropic conjugates over $K$ is a natural positive
complexity $c$ on the orbit of $z_0$ under
$\PSU_q(\OOO_K)$, which is invariant under $\PSU_q(\OOO_K)_\infty$.

\btheo\label{theo:countingcubicintro} There exists a constant
$\kappa>0$ and an explicit constant $C>0$ such that, as $s\ra+\infty$,
$$
\card\{z\in \PSU_q(\OOO_K)_\infty\bs 
\PSU_q(\OOO_K)\cdot z_0\;:\; c(z)\leq s\}=
C\;s^4\,(1+\bigO(s^{-\kappa}))\;.
$$
\etheo
We refer to Section \ref{sect:cubicpoints} for a more precise and more
general statement, valid for instance for the congruence subgroups of
$\PSU_q(\OOO_K)$.

\medskip The main tools of this paper, described in Section
\ref{sect:geometry}, are a geometric counting result for the common
perpendiculars between two properly immersed closed convex subsets,
and an equidistribution result of their initial and terminal tangent
vectors, valid in general pinched negative curvature
\cite{ParPau14}. See also \cite{OhSha12, OhShaCircles} for related
counting and equidistribution results in real hyperbolic spaces,
\cite{Kim13} in negatively curved symmetric spaces, \cite{ParPauHub}
for more geometric applications, and \cite{ParPauTou} for arithmetic
applications using real hyperbolic spaces. 

We recall the necessary geometric background on the complex hyperbolic
spaces in Section \ref{sec:cxhyp}. Section \ref{sect:complexhypcase}
is devoted to the computation of the measure theoretic constants that
appear in the geometric counting and equidistribution result in the
special case of the complex hyperbolic spaces, as was done in \cite[\S
7]{ParPauRev} in the case of the real hyperbolic spaces. In
particular, these computations allow to obtain the following geometric
counting result (see Corollary \ref{coro:complexhyperbo} for a more
general version).

\btheo Let $\Ga$ be a discrete group of isometries of the complex
hyperbolic $n$-space $\hnc$ such that the orbifold $M=\Ga\bs\hnc$ has
finite volume. Let $D^-$ be an horoball in $\hnc$ centred at a
parabolic fixed point of $\Ga$ and let $D^+$ be a complex geodesic
line in $\hnc$ whose quotient under its stabiliser $\Ga_{D^+}$ in
$\Ga$ has finite volume, and let $m^+$ be the cardinality of the
pointwise stabiliser of $D^+$ in $\Ga$. Then the number (counted with
multiplicities) of the common perpendiculars of length at most $t$
between the images of $D^-$ and $D^+$ in $M$ is equivalent as
$t\ra+\infty$ to
$$
\frac{4\,(n-1)}{(2n-1)\,\pi}
\frac{\Vol(\Ga_{D^-}\bs D^-)\Vol(\Ga_{D^+}\bs D^+)}
{m^+\;\Vol(M)}\;e^{2n\,t}\,.
$$
\etheo

In order to motivate our next results, let us state in an
appropriate way the classical result, known as {\it Mertens' formula},
describing the asymptotic behaviour of the average order of Euler's
function, or equivalently the asymptotic counting of Farey fractions,
and its related equidistribution result of Farey fractions in the
group $\RR$. The additive group $\ZZ$ acts on $\ZZ\times\ZZ$ by
horizontal shears (transvections): $k\cdot (u,v)=(u+kv,v)$. Then
$$
\card\;\;\;\ZZ\,\bs
\big\{(u,v)\in\ZZ\times\ZZ\;: (u,v)=1,\  |v|\le s \big\}=
\frac 6{\pi^2}\,s+\operatorname{O}(s^{1-\kappa})
$$
for some $\kappa>0$ (see for example \cite[Thm.~330]{HarWri08}, a
better error term is due to Walfisz \cite{Walfisz63}). Furthermore,
with $\Delta_{x}$ the unit Dirac mass at $x$, as $s\to+\infty$,
$$
\frac{\pi^2}{6\,s}\;\;\sum_{(u,\,v)=1,\;|v|\le s}\Delta_{\frac uv}\;
\weakstar\;\Leb_\RR\,.
$$

Our next result is an analog of Mertens' formula for Heisenberg
groups. Let $\tr,\n$ be the (absolute) trace and norm of $K$, and let
$\langle a,b,c\rangle$ be the ideal of $\OOO_K$ generated by $a,b,c\in
\OOO_K$.  The nilpotent group
$$
\Heis_3(\OOO_K)=\{(w_0,w)\in \OOO_K\times \OOO_K\;:\; \tr(w_0)=\n(w)\}
$$
with law \begin{equation}\label{eq:heislaw}
(w_0,w)(w'_0,w')= (w_0+w'_0+w'\,\overline{w},w+w')
\end{equation} acts on
$\OOO_K\times\OOO_K\times\OOO_K$ by the shears
$$
(w_0,w)(a,\alpha,c)=(a+\overline{w}\,\alpha+w_0\,c,\alpha+w\,c,c)\;.
$$

\btheo\label{theo:countHeisintro}  
There exists $\kappa>0$ such that, as $s\ra+\infty$,
\begin{align*}
&\card\;\;\Heis_3(\OOO_K)\bs
\big\{(a,\alpha,c)\in\OOO_K\times\OOO_K\times\OOO_K\;:\;\begin{array}{c}\langle
a,b,c\rangle=\OOO_K\,,\\
\tr(a\,\overline{c})=\n(\alpha), \;\n(c)\leq s\end{array}\big\}\\ =\; & 
\frac{\zeta(3)}
{2\,\pi\,|D_K|^{\frac{1}{2}}\,\zeta_K(3)}
\;s^2+\bigO(s^{2-\kappa})\,, 
\end{align*}
where $D_K$ is the discriminant and $\zeta_K$ Dedekind's zeta function
of $K$.  
\etheo

The $3$-dimensional Heisenberg group
$$
\Heis_{3}= \{(w_0,w)\in\CC\times\CC \;:\; 2\,\Re\;w_0
=|w|^2\}\,,
$$
with the group law 
\eqref{eq:heislaw} is
the Lie group of $\RR$-points of a $\QQ$-group whose group of
$\QQ$-points is $\Heis_{3}\cap (K\times K)$. We endow it with its Haar
measure 
\begin{equation}\label{eq:defhaarheis}
d\haarheis(w_0,w)= d(\Im\,w_0)\, d(\Re\,w)\, d(\Im\,w)\,.
\end{equation} 
Theorem \ref{theo:countHeisintro} is a counting result of  rational points
$(\frac{a}{c},\frac{\alpha}{c})$ (analogous to Farey fractions) in
$\Heis_3$, and the following result is a related equidistribution
theorem.

\btheo\label{theo:equidisHeisintro} As $s\ra+\infty$, we have
$$
\frac{\pi\,|D_K|^{\frac{3}{2}}\,\zeta_K(3)}
{\zeta(3)}\;s^{-2}\sum_{\substack{(a,\,\alpha,\,c)\in
\OOO_K\times\OOO_K\times\OOO_K,
\;0<\n(c)\leq s \\ \tr(a\,\overline{c})=\n(\alpha),
\;\langle a,\,\alpha,\,c\rangle=\OOO_K}}\;
\Delta_{(\frac{a}{c},\frac{\alpha}{c})}\;\weakstar\;\haarheis\,.
$$ 
\etheo

Theorems \ref{theo:countHeisintro} and \ref{theo:equidisHeisintro}
being of an arithmetic nature without any reference to the geometry
used in their proofs, can certainly be proven using techniques from
analytic number theory, and such a direct approach can produce
somewhat more precise results.  We refer to Theorems
\ref{theo:countHeis} and \ref{theo:equidisHeis} in Section
\ref{sect:mertensheisenberg} for more general results with added
congruence properties, and to Remark \ref{rem:highdim} for counting
and equidistribution results in higher dimensional Heisenberg groups.

\medskip\noindent{\small {\it Acknowledgement: } The first author
  thanks the Université de Paris-Sud (Orsay) for a visit of a month
  and a half which allowed an important part of the writing of this
  paper, under the financial support of the ERC grant GADA 208091. We
  thank Y.~Benoist and L.~Clozel for their help with Proposition
  \ref{prop:caracfixpoint}.}

\section{Geometric counting and equidistribution}
\label{sect:geometry}

In this Section, we briefly review a simplified version of the
geometric counting and equidistribution results proved in
\cite{ParPau14}, whose arithmetic applications will be considered in
the main part of this paper (see also \cite{ParPauRev} for a review of
related references).

Let $\wt M$ be a negatively curved rank one symmetric space (see
\cite[\S 2]{ParPau14} for a more general setting). In other words,
$\wt M$ is a hyperbolic space $\HH^n_\FF$ where $\FF$ is the set $\RR$
of real numbers, $\CC$ of complex numbers, $\HH$ of Hamilton's
quaternions, or $\OO$ of octonions, and $n\geq 2$, with $n=2$ if
$\KK=\OO$ see for instance \cite{Mostow73,Parker07}). We will
normalise them so that their maximal sectional curvature is $-1$.  

Let $\Ga$ be a discrete nonelementary group of isometries of $\wt M$ and
let $M=\Ga\bs\wt M$ be the quotient orbifold. Let $D^-$ and $D^+$ be
nonempty proper closed convex subsets of $\wt M$, such that the
families $(\ga D^-)_{\ga\in\Ga}$ and $(\ga D^+)_{\ga\in\Ga}$ are
locally finite in $\wt M$ (see \cite[\S 3.3]{ParPau14} for more
general families). Let $\Ga_{D^-}$ and $\Ga_{D^+}$ be the stabilisers
in $\Ga$ of the subsets $D^-$ and $D^+$, respectively.

\medskip We denote by $\partial_\infty \wt M$ the boundary at infinity
of $\wt M$, by $\Lambda \Ga$ the limit set of $\Ga$ and by
$(\xi,x,y)\mapsto \beta_\xi(x,y)$ the Busemann function on
$\partial_\infty \wt M\times \wt M\times \wt M$ (see for instance
\cite{BriHae99}). For every $v\in T^1\wt M$, let $\pi(v)\in \wt M$ be
its origin, and let $v_-, v_+$ be the points at infinity of the
geodesic line defined by $v$. We denote by $\normalpm D^\mp$ the {\it
  outer/inner unit normal bundle} of $\partial D^\mp$, that is, the
set of $v\in T^1\wt M$ such that $\pi(v)\in \partial D^\mp$ and the
closest point projection on $D^\mp$ of $v_\pm\in\partial_\infty \wt
M-\partial_\infty D^\mp$ is $\pi(v)$. For all $\ga,\ga'$ in $\Ga$, the
convex sets $\ga D^-$ and $\ga' D^+$ have a common perpendicular if
and only if their closures $\overline{\ga D^-}$ and $\overline{\ga'
  D^+}$ in $\wt M\cup\partial_\infty \wt M$ do not intersect.  We
denote by $\alpha_{\ga,\,\ga'}$ this common perpendicular (starting
from $\ga D^-$ at time $t=0$) and by $\ell(\alpha_{\ga,\,\ga'})$ its
length.  The {\em multiplicity} of $\alpha_{\ga,\ga'}$ is
$$
m_{\ga,\ga'}=
\frac 1{\card(\ga\Ga_{D^-}\ga^{-1}\cap\ga'\Ga_{D^+}{\ga'}^{-1})}\,,
$$
which equals $1$ when $\Ga$ acts freely on $T^1\wt M$ (for instance
when $\Ga$ is torsion-free). Let
$$
\N_{D^-,\,D^+}(t)=
\sum_{
(\ga,\,\ga')\in \Ga\bs((\Ga/\Ga_{D^-})\times (\Ga/\Ga_{D^+}))\;:\;
\overline{\ga D^-}\,\cap \,\overline{\ga' D^+}\,=\emptyset,\; 
\ell(\alpha_{\ga,\, \ga'})\leq t} m_{\ga,\ga'}
\;,
$$
where $\Ga$ acts diagonally on $\Ga\times\Ga$. When $\Ga$ has no
torsion, $\N_{D^-,\,D^+}(t)$ is the number (with multiplicities coming
from the fact that $\Ga_{D^\pm}\bs D^\pm$ is not assumed to be
embedded in $M$) of the common perpendiculars of length at most $t$
between the images of $D^-$ and $D^+$ in $M$. We refer to \cite[\S
4]{ParPau14} for the use of H\"older potentials on $T^1\wt M$ to
modify this counting function by adding weights, which could be useful
for some further arithmetic applications.

Recall the following notions (see for instance \cite{Roblin03}). The
{\em critical exponent} of $\Ga$ is
$$
\delta_{\Ga}=\limsup_{n\to+\infty}\frac 1n \ln\card\{\ga\in\Ga:
d(x_{0},\ga x_{0})\leq n\}\,,
$$
which is positive, finite and independent of $x_{0}\in\wt M$. Let
$(\mu_{x})_{x\in \wt M}$ be a {\em Patterson density} for $\Ga$, that
is a family $(\mu_{x})_{x\in \wt M}$ of nonzero finite measures on
$\partial_{\infty}\wt M$ whose support is $\Lambda\Ga$, such that
$\ga_*\mu_x=\mu_{\ga x}$ and
$$
\frac{d\mu_{x}}{d\mu_{y}}(\xi)=e^{-\delta_{\Ga}\beta_{\xi}(x,\,y)}
$$
for all $\ga\in\Ga$, $x,y\in\wt M$ and $\xi\in\partial_{\infty}\wt M$.

The {\em Bowen-Margulis measure} $\wt m_{\rm BM}$ for $\Ga$ on $T^1\wt
M$ is defined, using Hopf's parametrisation $v\mapsto
(v_-,v_+,\beta_{v_+}(x_0,\pi(v))\,)$ from $T^1\wt M$ into
$\partial_\infty \wt M\times\partial_\infty \wt M\times\RR$, by
$$
d\wt m_{\rm BM}(v)=e^{-\delta_{\Ga}(\beta_{v_{-}}(\pi(v),\,x_{0})+
\beta_{v_{+}}(\pi(v),\,x_{0}))}\;
d\mu_{x_{0}}(v_{-})\,d\mu_{x_{0}}(v_{+})\,dt\,.  
$$
Note that in the right hand side of this equation, $\pi(v)$ may be
replaced by any point $x'$ on the geodesic line defined by $v$, since
$\beta_{v_{-}}(\pi(v),\,x')+ \beta_{v_{+}}(\pi(v),\,x')=0$.  We will
use this elementary observation in the proof of Lemma
\ref{lem:computheisen}(ii).  The measure $\wt m_{\rm BM}$ is nonzero,
independent of $x_{0}\in\wt M$, is invariant under the geodesic flow,
the antipodal map $v\mapsto -v$ and the action of $\Ga$. Thus, it
defines a nonzero measure $m_{\rm BM}$ on $T^1M$ which is invariant
under the geodesic flow of $M$ and the antipodal map, called the {\em
  Bowen-Margulis measure} on $M=\Ga\bs \wt M$.  When $m_{\rm BM}$ is
finite (for instance when $M$ has finite volume or when $\Ga$ is
geometrically finite), denoting the total mass of a measure $m$ by
$\|m\|$, the probability measure $\frac{m_{\rm BM}} {\|m_{\rm BM}\|}$
is then uniquely defined, and is the unique probability measure of
maximal entropy for the geodesic flow (see \cite{OtaPei04}).

\medskip Using the endpoint homeomorphism $v\mapsto v_+$ from
$\normalout {D^-}$ to $\partial_{\infty}\wt M-\partial_{\infty}D^-$,
we defined in \cite{ParPau13ETDS} (generalising the definition of Oh
and Shah \cite[\S 1.2]{OhSha12} when $D^-$ is a horoball or a
totally geodesic subspace in $\hnr$) the {\em skinning measure}
$\wt\sigma_{D^-}$ of $\Ga$ on $\normalout{D^-}$, by
$$
d\wt\sigma_{D^-}(v) =  
e^{-\delta_\Ga\,\beta_{v_{+}}(\pi(v),\,x_{0})}\,d\mu_{x_{0}}(v_{+})\,.
$$
The measure $\wt\sigma_{D^-}$ is independent of $x_{0}\in\wt M$, is
nonzero if $\Lambda\Ga$ is not contained in $\partial_{\infty}D^-$,
and satisfies $\wt\sigma_{\ga D^-} =\ga_{*}\wt\sigma_{D^-}$ for every
$\ga\in\Ga$. Since the family $(\ga D^-)_{\ga\in\Ga}$ is locally
finite in $\wt M$, the measure $\sum_{\ga \in \Ga/\Ga_{D^-}}
\;\ga_*\wt\sigma_{D^-}$ is a well defined $\Ga$-invariant locally
finite (nonnegative Borel) measure on $T^1\wt M$.  Hence, it induces a
locally finite measure $\sigma_{D^-}$ on $T^1M=\Ga\bs T^1\wt M$,
called the {\em skinning measure} of $D^-$ in $T^1M$. We refer to
\cite[\S 5]{OhShaCircles} and \cite[Theo.~9]{ParPau13ETDS} for
finiteness criteria of the skinning measure $\sigma_{D^-}$, in
particular satisfied when $M$ has finite volume and if either $D^-$ is
a horoball centred at a parabolic fixed point of $\Ga$ or if $D^-$ is
a totally geodesic subspace.

\brema\label{lem:muinfty} {\rm Let $\H$ be a horoball in $\wt M$
  centred at $\xi$, and let $\rho$ be the geodesic ray starting from
  any point in $\partial \H$ and converging to $\xi$. Then (see for
  instance \cite[\S 2.3]{HerPau04}) the weak-star limit
$$
\mu_\xi=\lim_{t\ra+\infty}\; e^{\delta_{\Ga} t}\,\mu_{\rho(t)}
$$
exists, defines a measure on $\partial_{\infty}\wt M-\{\xi\}$ which is
invariant under the elements of $\Ga$ preserving $\H$. The limit
measure satisfies
\begin{equation}\label{eq:pattersoninfinity}
\frac{d\mu_x}{d\mu_\xi} (\eta)=e^{-\delta_{\Ga}\beta_{\eta}(x,\,x_{\H,\,\eta})}
\end{equation}
for all $x\in\wt M$, $\eta\in\partial_{\infty}\wt M-\{\xi\}$, where
$x_{\H,\,\eta}$ is the intersection with $\partial \H$ of the geodesic
line from $\eta$ to $\xi$.  Take $x_0=\rho(t)$ and let $t$ go to
$+\infty$ in the definition of the Bowen-Margulis measure and the
skinning measures. Then, for every $v\in T^1\wt M$ such that
$v_\pm\neq \xi$, we have
\begin{equation}\label{eq:BMinfinity}
d\wt m_{\rm BM}(v)=e^{-\delta_{\Ga}(\beta_{v_{-}}(\pi(v),\,x_{\H,\,v_-})+
\beta_{v_{+}}(\pi(v),\,x_{\H,\,v_+}))}\;
d\mu_{\xi}(v_{-})\,d\mu_{\xi}(v_{+})\,dt\,.
\end{equation}
Furthermore,  for every $v\in \normalout \H$, we have
\begin{equation}\label{eq:skinninginfinity}
d\wt\sigma_{\H}(v) = d\mu_{\xi}(v_{+})\,,
\end{equation}
since 
$\beta_{v_{+}}(\pi(v),\rho(t))=-t+\smallo(1)$ as $t\ra+\infty$.}
\erema

The following result is a special case of \cite[Coro.~20, 21,
Theo.~28]{ParPau14} (we prove that the pairs of the initial/terminal
tangent vectors of the common perpendiculars equidistribute in the
product of the outer/inner tangent bundles of $D^-$/$D^+$). We refer
to \cite{ParPauRev} for an review of the particular cases known before
\cite{ParPau14} (due to Huber, Margulis, Herrmann, Cosentino, Roblin,
Oh-Shah, Martin-McKee-Wambach, Pollicott, and the authors for
instance).

For all $t\geq 0$ and $x\in\partial D^-$, let 
$$
m_t(x)=\sum_{\ga\in \Ga/\Ga_{D^+}\;:\;
  \overline{D^-}\,\cap \,\overline{\ga D^+}\,= \emptyset,\;
\alpha_{e,\, \ga}(0)=x,\; \ell(\alpha_{e,\, \ga})\leq t} m_{e,\ga}
$$ 
be the multiplicity of $x$ as the origin of common perpendiculars with
length at most $t$ from $D^-$ to the elements of the $\Ga$-orbit of
$D^+$.  We denote by $\Delta_x$ the unit Dirac mass at a point $x$.

\btheo\label{theo:mainequicountdown} Let $\wt M,\Ga,D^-,D^+$ be as
above. Assume that the measures $m_{\rm BM},\sigma_{D^-},\sigma_{D^+}$
are nonzero and finite. Then
$$
\N_{D^-,\,D^+}(t)\;\sim\;
\frac{\|\sigma_{D^-}\|\;\|\sigma_{D^+}\|}{\delta_\Ga\;\|m_{\rm
    BM}\|}\; e^{\delta_\Ga t}\;,
$$
as $t\ra+\infty$. If $\Ga$ is arithmetic or if $M$ is compact, then
the error term is $\bigO(e^{(\delta_\Ga-\kappa) t})$ for some $\kappa
>0$.  Furthermore, the origins of the common perpendiculars
equidistribute in the boundary of $D^-$: 
\begin{equation}\label{eq:equidistribdown}
\lim_{t\ra+\infty}\; 
\frac{\delta_\Ga\;\|m_{\rm BM}\|}{\|\sigma_{D^-}\|\;\|\sigma_{D^+}\|}
\;e^{-\delta_\Ga t}\;\sum_{x\in\partial D^-} m_{t}(x) \;
\Delta_{x}\;=\; \frac{\pi_*\wt\sigma_{D^-}}{\|\sigma_{D^-}\|}\;
\end{equation}
for the weak-star convergence of measures on the locally compact space
$T^1\wt M$. 
\cqfd
\etheo

When $\partial D^-$ is smooth, for smooth functions $\psi$ with
compact support on $\partial D^-$, there is an error term in the
equidistribution claim \eqref{eq:equidistribdown} when the measures on
both sides are evaluated on $\psi$, of the form $\bigO(e^{-\kappa t}\,
\|\psi\|_\ell)$ where $\kappa>0$ and $\|\psi\|_\ell$ is the Sobolev
norm of $\psi$ for some $\ell\in\NN$, as proved in
\cite[Theo.~28]{ParPau14}.

\medskip When $M$ has finite volume, the Bowen-Margulis measure
$m_{\rm BM}$ coincides up to a multiplicative constant with the
Liouville measure on $T^1M$, and the skinning measures of points,
horoballs and totally geodesic subspaces $D^\pm$ coincide with the
(homogeneous) Riemannian measures on $\normalpm D^\mp$ induced by the
(Sasaki's) Riemannian metric of $T^1\wt M$. In Section
\ref{sect:complexhypcase}, we explicit some of these proportionality
constants when $\wt M$ is a complex hyperbolic space (see \cite[\S
7]{ParPauRev} for the real hyperbolic case).

\section{Complex hyperbolic geometry}
\label{sec:cxhyp}

In this Section, we recall some background on the complex hyperbolic
spaces, as mostly contained in \cite{Goldman99}, and, unless otherwise
stated, we will follow the conventions therein. For all $w,w'$ in
$\CC^{n-1}$, we denote by $w\cdot \overline{w'}= \sum_{i=1}^{n-1}
w_i\,\overline{w'_i}$ their standard Hermitian product, and we denote
$|w|^2= w\cdot \overline w$.  Recall that for every $n\geq 1$, the
Siegel domain model of the complex hyperbolic $n$-space $\hnc$ is
$$ 
\big\{(w_0,w)\in\CC\times\CC^{n-1}\;:\; 
2\,\Re\; w_0 -|w|^2>0\big\}\,,
$$ 
endowed with the Riemannian metric
$$
ds^2_{\,\hnc}=\frac{1}{(2\,\Re\; w_0 -|w|^2)^2}
\big((dw_0-dw\cdot \overline{w})((\overline{dw_0}-w\cdot 
\overline{dw})+
(2\,\Re\; w_0 -|w|^2)\;dw\cdot \overline{dw}\big)\,.
$$
In accordance with Section \ref{sect:geometry}, this metric is
normalised so that its sectional curvatures are in $[-4,-1]$, instead
of in $[-1,-\frac{1}{4}]$ as in \cite{Goldman99} and \cite{Parker12}.
Its boundary at infinity is
$$
\partial_\infty\hnc=\big\{(w_0,w)\in
\CC\times \CC^{n-1} \;:\; 2\,\Re\; w_0 -|w|^2=0\big\}\cup\{\infty\}\,.
$$
A {\it complex geodesic line} in $\hnc$ is the image by an isometry of
$\hnc$ of the intersection of $\hnc$ with the complex line
$\CC\times\{0\}$; with our normalisation of the metric, a complex line
has constant sectional curvature $-4$. The boundary at infinity of a
complex geodesic line is a topological circle, called a {\it chain}
(see Section \ref{sect:chains} for more informations). A chain is {\em
  finite} if it does not contain $\infty$.

\medskip
Let $q$ be the nondegenerate Hermitian form $-z_0\,\overline{z_n}
-z_n\,\overline{z_0} + |z|^2$ of signature $(1,n)$ on
$\CC\times\CC^{n-1}\times\CC$ with coordinates $(z_0,z,z_n)$. This is
not the form considered in \cite[p.~67]{Goldman99}, hence we need to
do some computations with it, but it is better suited for our
purposes. It is the one considered for instance in \cite{ParPau10GT},
to which we will refer frequently.  The Siegel domain $\hnc$ embeds in
the complex projective $n$-space $\PP_n(\CC)$ by the map (using
homogeneous coordinates)
$$
(w_0,w)\mapsto [w_0:w:1]\;.
$$
We identify $\hnc$ with its image by this map. This image, called the
{\it projective model} of $\hnc$ when endowed with the isometric
Riemannian metric, is the {\it negative cone} of $q$, that is
$\{[z_0:z:z_n]\in \PP_n(\CC) \;:\;q(z_0,z,z_n)<0\}$. This embedding
extends continuously to the boundary at infinity, by mapping
$(w_0,w)\in
\partial_\infty\hnc-\{\infty\}$ to $[w_0:w:1]$ and $\infty$ to
$[1:0:0]$, so that the image of $\partial_\infty\hnc$ is the {\it null
  cone} of $q$, that is $\{[z_0:z:z_n]\in\PP_n(\CC)\;:\;
q(z_0,z,z_n)=0\}$. 

The linear action of the special unitary group of $q$
$$
\SU_q=\{g\in \SL_{n+1}(\CC)\;:\; q\circ g=q\}
$$
on $\CC^{n+1}$ induces a projective action on $\PP_n(\CC)$. The
quotient group $\PSU_q= \SU_q/(\UU_{n+1}\operatorname{Id})$ of $\SU_q$
by the kernel of the projective action, where $\UU_{n+1}$ is the group
of $(n+1)$-th roots of unity, preserves $\hnc$, and its restriction to
$\hnc$ is the orientation-preserving isometry group of $\hnc$. For
instance by the paragraph above \cite[Lem.~6.3]{ParPau10GT}, an
element $\ga\in\SU_q$ fixes $\infty$ if and only if $\ga$ is upper
triangular (this is the reason, besides rationality problems, that we
chose the Hermitian form $q$ rather than the one in \cite{Goldman99}),
see for instance \cite[p.~119]{Goldman99}, \cite[\S 2.1]{FalPar06} up
to signs.

By for instance \cite[Eq.~(42)]{ParPau10GT}, the intersection of
$\SU_q$ with the upper triangular subgroup of $\SL_3(\CC)$ is 
$$
B_q=\left\{\left(\!\!\begin{array}{ccc}
a_1&\overline{\zeta}&\frac{1}{2\,\overline{a_1}}(|\zeta|^2-iu)\\
0 & a_2 & \frac{a_2}{\overline{a_1}}\,\zeta\\
0 & 0 & a_3
\end{array}\!\!\right)\;:\;\begin{array}{c}
\zeta,a_1,a_2,a_3\in\CC,\;u\in\RR,\\
a_1a_2a_3=1, \,a_3\,\overline{a_1}=1,\;|a_2|=1\end{array}\right\}\,,
$$
and its image $\overline{B_q}$ in $\PSU_q$ is equal to the stabiliser
in $\PSU_q$ of $\infty$. 

If $\ga\in \SU_q$ induces a loxodromic isometry on $\hnc$, then $\ga$
is diagonalisable over $\CC$, it has a unique eigenvalue $\lambda$ of
modulus $>1$, and its eigenvalues are $\lambda, \frac{1}
{\overline{\lambda}}, \frac{\overline{\lambda}}{\lambda}$ (see for
instance \cite[\S 3.2]{Parker12}). Furthermore, the translation length
$\ell$ of $\ga$ in $\hnc$ is
\begin{equation}\label{eq:calctransllengthloxo}
\ell= \ln |\lambda|\,,
\end{equation}
(see for instance \cite[Prop.~3.10]{Parker12}, noting that this
reference normalises the curvature to be between $-1$ and $-1/4$).

\medskip
The {\it horospherical coordinates} $(\zeta,u,t)\in
\CC^{n-1}\times\RR\times\,[0,+\infty[$, that we will use from now on
unless otherwise stated, of $(w_0,w)\in\hnc\cup (\partial_\infty \hnc-
\{\infty\})$ are
\begin{equation}\label{eq:horosphecoord}
(\zeta,u,t)=(w,\;-2\,\Im\; w_0, \;2\,\Re\; w_0-|w|^2)
\;\;\;{\rm hence}\;\;\;
(w_0,w)=\big(\frac{1}{2}(|\zeta|^2+t-iu) ,\zeta\big)\,,
\end{equation}
so that the Riemannian metric is given by 
\begin{equation}\label{eq:riemethorosphecoord}
ds^2_{\,\hnc}= \frac{1}{4\,t^2}\big(dt^2
+(du+2\,\Im\;d\zeta\cdot\overline{\zeta}\,)^2+
4\,t\,d\zeta\cdot\overline{d\zeta} \,\big)\,.
\end{equation}  
In horospherical coordinates, the geodesic lines from
$(\zeta,u,0)\in\partial_\infty \hnc- \{\infty\} $ to $\infty$ are, up
to translations at the source, the map $s\mapsto (\zeta,u,e^{2s})$, by
the normalisation of the metric.  The closed horoballs centred at
$\infty\in\partial_\infty \hnc$ are the subsets
$$
\H_s=\{(\zeta,u,t)\in\hnc\;:\;t\geq s\},
$$
and the horospheres centred at $\infty$ are their boundaries
$$
\partial\H_s=\{(\zeta,v,t)\in\hnc\;:\;t=s\}\,,
$$ 
for any $s>0$. Note that, for every $s\geq 1$, we have
\begin{equation}\label{eq:distentrhorob}
d(\partial\H_1,\partial\H_s)=\frac{\ln s}{2}\;.
\end{equation}

As introduced by \cite[p.~297]{Parker92}, the {\it Cygan distance} on
$\hnc\cup (\partial_\infty \hnc-\{\infty\})$ (analogous to the
Euclidean distance on the closure in $\RR^n$ of the upper halfspace
model of $\hnr$) is
$$
d_{\rm Cyg}((\zeta,u,t),(\zeta',u',t'))=
\big|\,|\zeta-\zeta'|^2+|t-t'|+
i(u-u'+2\,\Im\;\zeta\cdot\overline{\zeta'})\big|^{1/2}\,.
$$

The {\it Heisenberg group} $\Heis_{2n-1}$ of dimension $2n-1$ is
the real Lie group structure on $\CC^{n-1}\times \RR$ with law
$$
(\zeta,u)(\zeta',u')=
(\zeta+\zeta',u+u'+2\,\Im\;\zeta\cdot\overline{\zeta'})
$$
and inverses $(\zeta,u)^{-1}=(-\zeta,-u)$. It identifies with
$\partial_\infty \hnc-\{\infty\}$ by the map $(\zeta,u)\mapsto
(\zeta,u,0)$. It acts on $\hnc\cup (\partial_\infty \hnc-\{\infty\})$
by the {\it Heisenberg translations}
$$
(\zeta,u)(\zeta',u',t')= (\zeta+\zeta',u+u'+ 2\,\Im\;
\zeta\cdot\overline{\zeta'},t')\;,
$$
that are isometries for both the Riemannian metric and the Cygan
distance, and that preserve the horospheres centred at $\infty$.  For
every $u\in\RR$, the Heisenberg translation by $(0,u)$ is called a
{\it vertical translation}.

It is easy to see that the Cygan distance\footnote{It is called the
  {\it Kor\'anyi distance} by many people working in sub-Riemannian
  geometry, though Kor\'anyi \cite{Koranyi85} does attribute it to
  Cygan \cite{Cygan78}.} on $\Heis_{2n-1}$ (see \cite[page
160]{Goldman99}) is the unique left-invariant distance on
$\Heis_{2n-1}$ with $d_{\rm Cyg} ((\zeta,u),(0,0))= (|\zeta|^4
+u^2)^{\frac{1}{4}}$.  We introduced in \cite[\S 6.1]{ParPau10GT} the
{\it modified Cygan distance} $d'_{\rm Cyg}$ as the unique
left-invariant distance on $\Heis_{2n-1}$ with
$$
d'_{\rm Cyg}((\zeta,u),(0,0))=
((|\zeta|^4+u^2)^{\frac{1}{2}}+|\zeta|^2)^{\frac{1}{2}}\,.
$$
We introduced in \cite[Lem.~3.4]{ParPau11MZ} the map $d''_{\rm Cyg}=
\frac{{d_{\rm Cyg}}^2}{d'_{\rm Cyg}}$, which is almost a distance on
$\Heis_{2n-1}$, as it satisfies $\frac{1}{\sqrt{2}}\;d_{\rm Cyg}\leq
d''_{\rm Cyg} \leq d_{\rm Cyg}$. For every nonempty subset $A$ of
$\Heis_{2n-1}$, we define the diameter of $A$ for this almost distance
as
$$
\diam_{d''_{\rm Cyg}}(A)=\sup_{x,\,y\,\in A} d''_{\rm Cyg} (x,y)\;.
$$

The following result on diameters of finite chains will be useful in
Section \ref{sect:chains}.

\blemm\label{lem:diamchain} For every finite chain $C$, we have
$$
\diam_{d_{\rm Cyg}}(C)=\frac{1}{\sqrt{2}}\;\diam_{d'_{\rm Cyg}}(C) =
\sqrt{2}\;\diam_{d''_{\rm Cyg}}(C)\;.
$$
\elemm

\dem For every chain $C$, there exists (see for instance \cite[\S
3.1.4]{Goldman99}) a unique point $P=[z_0:z:z_n]$ in $\PP_n(\CC)$,
called the {\it polar point of} $C$, such that $q(z_0:z:z_n)>0$ and
$C$ is the intersection with $\partial_\infty\hnc$ of the orthogonal
to $P$ for $q$. Note that if $\ga\in\PSU_q$, then the polar point of
$\ga C$ is $\ga P$, and that $C$ is finite is and only if $z_n\neq 0$.

Let $P=[z_0:z:z_n]$ be the polar point of a finite chain $C$. Let
$\ga\in\PSU_q$ be the Heisenberg translation by
$$ 
\Big[\;\frac{|z|^2}{2\,|z_n|^2} +i\,\Im\;\frac{z_0}{z_n}\,:
\,-\frac{z}{z_n}\,:\,1\Big]\in \Heis_{2n-1}\;.
$$
An easy computation shows that $\ga[z_0:z:z_n]=[-a:0:b]$ with
$a=q(z_0,z,z_n)>0$ and $b=2|z_2|^2>0$. Since the distances $d_{\rm
  Cyg}$, $d'_{\rm Cyg}$ and $d''_{\rm Cyg}$ are left-invariant, we may
assume that $P$ is the point $[-R^2/2:0:1]$ with $R>0$. Hence, using
the facts that $\Re\;w_0=|w|^2$, $w=\zeta$ and
$\Im\;w_0=-\frac{u}{2}$, we have
$$
C=\{[w_0:w:1]\in\Heis_{2n-1}\,:\; R^2/2 -\omega_0=0\}
=\{(\zeta,u)\in\Heis_{2n-1}\,:\; u=0,\;|\zeta|=R\}\;,
$$
which is the sphere with radius $R$ in the first factor $\CC^{n-1}$ of
$\Heis_{2n-1}$. Since the Heisenberg dilations $(\zeta,u)\mapsto
(\lambda\zeta,\lambda^2 u)$ with $\lambda >0$ are homotheties of ratio
$\lambda$ for $d_{\rm Cyg}$, $d'_{\rm Cyg}$ and $d''_{\rm Cyg}$, we
may assume that $R=1$. This proves in particular that the three
diameters are proportional, and we now compute the proportionality
constants.

For every $(\zeta,0)\in C$, we have $d_{\rm Cyg}((\zeta,0),(0,0))=1$
and $d'_{\rm Cyg}((\zeta,0),(0,0))=\sqrt{2}$, hence by the triangle
inequality, $\diam_{d_{\rm Cyg}}(C)\leq 2$ and $\diam_{d'_{\rm Cyg}}
(C)\leq 2\sqrt{2}$. If $\zeta=(1,0,\dots,0)$, we have
$$
d_{\rm Cyg}((\zeta,0),(-\zeta,0))=d_{\rm Cyg}((2\zeta,0),(0,0))=2
$$ 
and similarly $d'_{\rm Cyg}((\zeta,0),(-\zeta,0))=2\sqrt{2}$. Hence
$\diam_{d_{\rm Cyg}}(C)=2$ and $\diam_{d'_{\rm Cyg}} (C)= 2\sqrt{2}$,
which proves the first equality.

Assume first that $n\ge 3$. Since the group $U(n-1)$ acts transitively
on the unit sphere of $\CC^{n-1}$, since the stabiliser of
$(1,0,\dots,0)$ in $U(n-1)$ acts transitively on the complex planes
containing $(1,0,\dots,0)$ and on the real halflines of the second
factor of $\CC^{n-1}$, for all $\zeta,\zeta'\in\SSS^{2n-3}$, there
exist $\ga\in U(n-1)$, $\theta\in[0,2\pi]$ and $\varphi\in [0,\pi]$
such that
$$
(\ga\zeta,\ga\zeta')= 
\big((1,0,\dots,0),(e^{i\theta}\cos\varphi,\sin\varphi,0,\dots,0)\big)\;.
$$
Hence we may assume that $n=3$, and with $\zeta_{\theta,\varphi}=
(e^{i\theta}\cos\varphi,\sin\varphi)$, we have
$$
\diam_{d''_{\rm Cyg}}(C)=\min_{\theta\in[0,2\pi],\;\varphi\in [0,\pi]}
d''_{\rm Cyg}((\zeta_{0,0},0),(\zeta_{\theta,\varphi},0))\;.
$$
Note that 
$$
(-\zeta_{0,0},0)\cdot(\zeta_{\theta,\varphi},0)=
(\zeta_{\theta,\varphi}-\zeta_{0,0},
-2\,\Im\;\zeta_{0,0}\cdot\overline{\zeta_{\theta,\varphi}})
=((e^{i\theta}\cos\varphi-1,\sin\varphi),2\sin\theta\cos\varphi)\;,
$$
and that $|(e^{i\theta}\cos\varphi-1,\sin\varphi)|^2=
2(1-\cos\theta\cos\varphi)$. Therefore
\begin{align*}
d''_{\rm Cyg}((\zeta_{0,0},0),(\zeta_{\theta,\varphi},0))^2&=
d''_{\rm Cyg}((-\zeta_{0,0},0)\cdot(\zeta_{\theta,\varphi},0), (0,0))^2\\&=
\frac{4(1-\cos\theta\cos\varphi)^2+4\sin^2\theta\cos^2\varphi}
{(4(1-\cos\theta\cos\varphi)^2+4\sin^2\theta\cos^2\varphi)^{1/2}+
2(1-\cos\theta\cos\varphi)}\\&=
\frac{2}{\frac{1}{(1+\cos^2\varphi-2\cos\theta\cos\varphi)^{1/2}}+
\frac{1-\cos\theta\cos\varphi}
{1+\cos^2\varphi-2\cos\theta\cos\varphi}}\;.
\end{align*}
Now $1+\cos^2\varphi-2\cos\theta\cos\varphi\leq
2(1-\cos\theta\cos\varphi)\leq 4$ and both equalities hold for
instance if $\theta=\pi$ and $\varphi=0$. Hence $d''_{\rm Cyg}
((\zeta_{0,0},0),(\zeta_{\theta,\varphi},0))^2\leq 2$ with equality if
$\theta=\pi$ and $\varphi=0$, therefore $\diam_{d''_{\rm Cyg}}(C)
=\sqrt{2}$. This proves the result when $n\ge 3$. 

If $n=2$, for all $\zeta,\zeta'\in\SSS^{1}$, there exist $\ga\in
U(1)$ and $\theta\in[0,2\pi]$ such that $(\ga\zeta,\ga\zeta')=
(1,e^{i\theta})$ and the result follows from the same computations as
above with $\varphi=0$.  \cqfd

\medskip We conclude this section by two geometric lemmas that will be
useful in Section \ref{sect:complexhypcase}. See also \cite[\S
3]{Kim13}, with slightly different conventions, for a computation
similar to Lemma \ref{lem:busemann} based on \cite[p.~113]{Goldman99}. 
The Cygan distance, the
Poisson kernel $e^{\beta_{(\xi,\,r)}}$, the Patterson measures $\mu_x$
computed in the next section, and related quantities are useful in
potential theory on the Heisenberg group and for the study of the
hypoelliptic Laplacian in sub-Riemannian geometry, see for instance
\cite{FolSte74,Krantz09}.

\blemm\label{lem:busemann}
For all $x=(\zeta,u,t)$ and $x'=(\zeta',u',t')$ in $\hnc$, for all
$(\xi,r)\in \Heis_{2n-1}=\partial_\infty \hnc-\{\infty\}$, we have
$$
\beta_{(\xi,\,r)}(x,x')=\frac{1}{2}\,\ln 
\frac{t'\;d_{\rm Cyg}(x,(\xi,r))^4}{t\;d_{\rm Cyg}(x',(\xi,r))^4}\,.
$$
\elemm

\dem 
Since the Busemann function $\beta_{\infty}(x,x')$ is unchanged if
$x$, $x'$ are replaced by other points on the horospheres centred at
$\infty$ through them, and since the map $s\mapsto (0,0,e^{2s})$ is a
geodesic line in $\hnc$ from $(0,0)\in \partial_\infty\hnc$ to
$\infty$, we have 
$$
\beta_\infty(x,x')=\frac{1}{2}\,\ln \frac{t'}{t}\;.
$$

It is easy to check that the map $\iota:(w_0,w)\mapsto (\frac{1}{w_0},
\frac{w}{w_0})$ is an isometric involution of $\hnc$ sending $(0,0)
\in \partial_\infty\hnc$ to $\infty$. Hence, with $x=(w_0,w)$ and
$x'=(w'_0,w')$, using Equation \eqref{eq:horosphecoord} and the fact
that $d_{\rm Cyg}(x,(0,0))^4=4|w_0|^2$ and $d_{\rm Cyg} (x',(0,0))^4 =
4|w_0'|^2$, we have
$$
\beta_{(0,\,0)}(x,x')=\beta_{\iota(0,\,0)}(\iota x,\iota x')=
\frac{1}{2}\,\ln \frac{2\,\Re\; \frac{1}{w'_0}-|\frac{w'}{w'_0}|^2}
{2\,\Re\; \frac{1}{w_0}-|\frac{w}{w_0}|^2}
=\frac{1}{2}\,
\ln\;\frac{t'\;d_{\rm Cyg}(x,(0,0))^4}{t\;d_{\rm Cyg}(x',(0,0))^4}\;.
$$
The Heisenberg translation $\tau$ by $(\xi,r)$ preserves the last
horospherical coordinates and the Cygan distances. Thus,
$\beta_{(\xi,\,r)}(x,x')= \beta_{(0,\,0)} (\tau^{-1}x, \tau^{-1}x')$,
which implies the claim.  
\cqfd

\blemm \label{lem:orthprojgeod}
The orthogonal projection from $\partial_\infty\hnc-\{(0,0),
\infty\}$ to the geodesic line in $\hnc$ with points at infinity
$(0,0)$ and $\infty$ is $(w_0,w)\mapsto (2\,|w_0|,0)$, that is,
in horospherical coordinates, $(\zeta,u,0)\mapsto
(0,0,(|\zeta|^4+{u}^2)^{1/2})$.  
\elemm

In particular, the preimages by this orthogonal projection are the
spheres of center $(0,0)$ for the Cygan distance on
$\Heis_{2n-1}$. They are spinal spheres with complex spine
$\{(w_0,w)\in\hnc\;:\;w=0\}$, see \cite[\S 5.1.9]{Goldman99}.

\medskip \dem 
For every parameter $a$ ranging in $]0,+\infty[\,$, consider the
horosphere $\partial\H_a$ centred at $\infty$. Its image by the
isometric involution $\iota:(w,w_0)\mapsto (\frac{1}{w_0},
\frac{w}{w_0})$ is, using Equation \eqref{eq:horosphecoord}, the
horosphere $\{(\xi,r,t)\in\hnc\;:\;t=\frac{a}4((|\xi|^2+t)^2+r^2)\}$
centred at $(0,0)$. The image of this horosphere by the Heisenberg
translation by $(\zeta,u)$ is the horosphere
$$
\{(\xi,r,t)\in\hnc\;:\;t=\frac{a}4\,((|-\zeta+\xi|^2+t)^2+(-u+r-2\,\Im
\;\zeta\cdot\overline{\xi})^2)\}
$$ 

\noindent
\begin{minipage}{9.9cm}centred at $(\zeta,u)$. 
The orthogonal projection of $(\zeta,u)$ on
the geodesic line $\ell$ from $(0,0)$ to $\infty$ is attained when the
parameter $a$ gives a double point of intersection $(0,0,t)$ between
this horosphere and $\ell$. The quadratic equation 
$$
t=\frac{a}4\,((|\zeta|^2+t)^2+u^2)
$$ 
with unknown $t$ has a double solution if
and only if its reduced discriminant $\Delta'=
(|\zeta|^2-\frac{2}{a})^2- (|\zeta|^4 +u^2)$ vanishes, that is, since
$a>0$, if and only if $a=\frac{2}{(|\zeta|^4 +u^2)^{1/2} +
  |\zeta|^2}$, giving $t=(|\zeta|^4+u^2)^{1/2}$.  The result follows. \cqfd
\end{minipage}
\begin{minipage}{5cm}
\begin{center}
\begin{picture}(0,0)%
\includegraphics{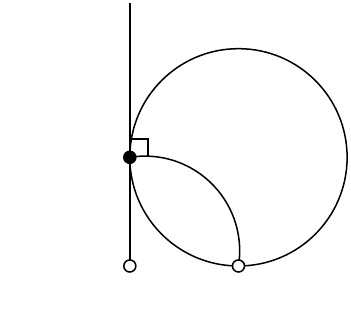}%
\end{picture}%
\setlength{\unitlength}{3812sp}%
\begingroup\makeatletter\ifx\SetFigFont\undefined%
\gdef\SetFigFont#1#2#3#4#5{%
  \reset@font\fontsize{#1}{#2pt}%
  \fontfamily{#3}\fontseries{#4}\fontshape{#5}%
  \selectfont}%
\fi\endgroup%
\begin{picture}(1733,1566)(1606,-1750)
\put(1621,-1006){\makebox(0,0)[lb]{\smash{{\SetFigFont{11}{13.2}{\rmdefault}{\mddefault}{\updefault}{\color[rgb]{0,0,0}$(0,0,t)$}%
}}}}
\put(2026,-1681){\makebox(0,0)[lb]{\smash{{\SetFigFont{11}{13.2}{\rmdefault}{\mddefault}{\updefault}{\color[rgb]{0,0,0}$(0,0)$}%
}}}}
\put(2566,-1681){\makebox(0,0)[lb]{\smash{{\SetFigFont{11}{13.2}{\rmdefault}{\mddefault}{\updefault}{\color[rgb]{0,0,0}$(\zeta,u)$}%
}}}}
\put(2116,-376){\makebox(0,0)[lb]{\smash{{\SetFigFont{11}{13.2}{\rmdefault}{\mddefault}{\updefault}{\color[rgb]{0,0,0}$\ell$}%
}}}}
\end{picture}%

\end{center}
\end{minipage}

\section{Measure computations in complex hyperbolic 
spaces}
\label{sect:complexhypcase}

In this Section, we give proportionality constants relating, on the
one hand, Patterson, Bowen-Margulis and skinning measures associated
to some convex subsets and, on the other hand, the corresponding
Riemannian measures, in the complex hyperbolic case.

We will denote the standard Lebesgue measure on $\CC^{n-1}$ by
$d\zeta$, and the usual left Haar measure on $\Heis_{2n-1}$ by
$$
d\lambda_{2n-1}(\zeta,u) = d\zeta du\,.
$$ 
In horospherical coordinates, the volume form of $\hnc= \Heis_{2n-1}
\times \;\mathopen{]}0,+\infty[$ is
\begin{equation}\label{eq:volhnc}
d\operatorname{vol}_{\hnc}(\zeta,u,t)=
\frac{1}{4\,t^{n+1}}\;d\zeta\,du\,dt\;.
\end{equation}

We begin by recalling a lemma that relates the volume of a Margulis
cusp neighbourhood with the volume of its boundary (compare \cite[Lem.~5.2]{HerPau96}).

\blemm\label{lem:horobvol}
Let $D$ be a horoball in $\hnc$ and let $\Ga$ be a discrete group of
isometries of $\hnc$ preserving $\partial D$ (hence $D$). Then
$\Vol(\Ga\backslash\partial D)=2n\Vol(\Ga\backslash D)$.  
\elemm 

\dem Since the group of isometries of $\hnc$ acts transitively on the
set of horospheres of $\hnc$, we may assume that $D=\H_1$.  The
horosphere centred at $\infty$ through a point $(\zeta,u,t) \in \H_1$
is $\partial \H_t$ and its orthogonal geodesic line at this point is
$s\mapsto (\zeta,u,e^{2s})$, hence
$$
d\vol_{\hnc}(\zeta,u,t)=d\vol_{\partial \H_t} (\zeta,u,t)
\,\frac{dt}{2t}\,.
$$  
By Equation \eqref{eq:volhnc}, we hence have 
\begin{equation}\label{eq:volhorosphzetau}
d\vol_{\partial \H_t}(\zeta,u,t)=\frac{1}{2\,t^{n}}\;d\zeta\,du
\end{equation} 
for every $t>0$, therefore $d\vol_{\partial \H_t}
(\zeta,u,t)=\frac{1}{t^{n}}\; d\vol_{\partial \H_1} (\zeta,u,1)$. The
homeomorphism from $\partial \H_t$ to $\partial \H_1$ defined by
$(\zeta,u,t)\mapsto (\zeta,u,1)$ commutes with the action of $\Ga$.
Thus,
\begin{align*}
\Vol(\Ga\backslash\H_1) &=\int_{\Ga\backslash\H_1}
\dvol_{\hnc}(\zeta,u,t)=\int_{t=1}^{+\infty}\int_{\Ga\backslash\partial \H_t}
\dvol_{\partial \H_t} (\zeta,u,t)\,\frac{dt}{2t}
\\ & =\int_{t=1}^{+\infty}\int_{\Ga\backslash\partial \H_1}
\dvol_{\partial \H_1} (\zeta,u,1)\,\frac{dt}{2t^{n+1}}=
\frac{1}{2n}\,\Vol(\Ga\backslash \partial \H_1)
\;. \;\;\;\Box
\end{align*}

\medskip
Let $\Ga$ be a lattice in $\Isom(\hnc)$. Its critical exponent is
$\delta_\Ga=2n$ (see for instance \cite[\S 6]{CorIoz99}). The
Patterson density $(\mu_x)_{x\in\hnc}$ of $\Ga$ is uniquely defined up
to a multiplicative constant, and is independent of $\Ga$. We will
choose the normalisation such that the measure $\mu_\infty$ defined 
in Remark \ref{lem:muinfty} by the horoball $\H_1$ 
coincides with $\lambda_{2n-1}$ on $\partial_\infty \hnc-\{\infty\}=
\Heis_{2n-1}$, that is 
$$
d\mu_\infty(\xi,r)= d\lambda_{2n-1}(\xi,r) =d\xi\,dr\;.
$$
This is possible since $\mu_\infty$ is invariant under the isometry
group $\Heis_{2n-1}$ preserving $\H_1$, hence it is a Haar measure on
$\partial_\infty \hnc-\{\infty\} = \Heis_{2n-1}$.

\blemm \label{lem:computheisen} Let $\Ga$ be a lattice in
$\Isom(\hnc)$, and let $\mu_\infty$ be normalised as above.  For all
$x=(\zeta,u,t)$ and $x'=(\zeta',u',t')$ in $\hnc$, for all
$(\xi,r)\in \partial_\infty \hnc-\{\infty\}$ and $v\in T^1\hnc$ such
that $v_\pm\neq\infty$, we have

\smallskip\noindent
(i)~ 
$\displaystyle d\mu_x(\xi,r)=
\frac{t^{n}}{d_{\rm Cyg}(x,(\xi,r))^{4n}}\;d\xi\,dr\;;$

\medskip\noindent
(ii) using a Hopf parametrisation $v\mapsto (v_-,v_+,s)$, 
$$
d\wt m_{\rm BM}(v)=
\frac{d\lambda_{2n-1}(v_-)\,d\lambda_{2n-1}(v_+)\,ds} 
{d_{\rm  Cyg}(v_-,v_+)^{4n}}\,;
$$

\noindent(iii) 
$$
\wt m_{\rm BM}= \frac{2n-1}{2^{2n-2}}\; \vol_{T^1\hnc}\,,
$$
and in particular 
$$
\|m_{\rm BM}\|= \frac{(2n-1)\,\pi^{n}}{2^{2n-3}\,(n-1)!}\;\Vol(M)\;;
$$

\medskip\noindent(iv) using the homeomorphism $v\mapsto v_+$ from
$\normalout\H_1$ to $\partial_\infty \hnc-\{\infty\}= \Heis_{2n-1}$,
$$
d\wt\sigma_{\H_1}(v)=d\lambda_{2n+1}(v_+)\,,
$$
and for every horoball $D^-$ in $\hnc$, 
$$
\pi_*\wt\sigma_{D^-}=2\,\vol_{\partial D^-}\,,
$$
and
$$
\|\sigma_{D^-}\|= 4n\, \Vol(\Ga_{D^-}\bs D^-)\;;
$$

\medskip\noindent (v) for every geodesic line $D^-$ in $\hnc$,
$$
d\pi_*\wt \sigma_{D^-}=\frac{n}{4^{n-1}\,(2n-1)}\;
d\pi_*\vol_{\normalout D^-}
$$
and, with $m$ the order of the pointwise stabiliser of $D^-$ in $\Ga$,
$$
\|\sigma_{D^-}\|= 
\frac{2\,\pi^{n-1}\, n!}{m\,(2n-1)!}\;\Vol(\Ga_{D^-}\bs D^-)\,;
$$

\medskip\noindent (vi) for every complex geodesic line $D^-$ in
$\hnc$, we have
$$
d\pi_*\wt \sigma_{D^-}= \frac{1}{2^{2n-3}}\;d\pi_*\vol_{\normalout D^-}
$$
and, with $m$ the order of the pointwise stabiliser of $D^-$ in $\Ga$,
$$
\|\sigma_{D^-}\|= 
\frac{\pi^{n-1}}{m\,4^{n-2}\,(n-2)!}\;\Vol(\Ga_{D^-}\bs D^-)\,.
$$  
\elemm

\dem In the computations below, it is useful to note that Lemma
\ref{lem:busemann} implies that
\begin{equation}\label{eq:rewritelem7}
e^{-2n\,\beta_{(\xi,\,r)}(x,\,x')}= \frac{t^n\;d_{\rm Cyg}(x',(\xi,r))^{4n}}
{(t')^n\;d_{\rm Cyg}(x,(\xi,r))^{4n}}\,.
\end{equation}

\medskip\noindent {\em (i)~} The geodesic line from $(\xi,r)$ to
$\infty$ goes through $\partial \H_1$ at the point $(\xi,r,1)$. By the
normalisation of $d\mu_\infty$ and by Equation
\eqref{eq:pattersoninfinity}, we hence have
$$
\frac{d\mu_x}{d\xi\,dr} (\xi,r)=
\frac{d\mu_x}{d\mu_\infty}(\xi,r)
=e^{-2n\,\beta_{(\xi,\,r)}(x,\,(\xi,\,r,\,1))}\,.
$$
The result then follows from Equation \eqref{eq:rewritelem7}.

\medskip\noindent {\em (ii)~} Note that if $x'$ is on the geodesic
line defined by $v$, then $d_{\rm Cyg} (x',v_-)^2\sim t'$ as $x'\ra
v_-$.  Hence, by Equation \eqref{eq:rewritelem7} and Assertion {\em
  (i)}, by letting $x'$ converge to $v_-$ on the geodesic line defined
by $v$, we have
\begin{align*}
&d\wt m_{\rm BM}(v) 
=e^{-2n(\beta_{v_{-}}(x',\,x)+ \beta_{v_{+}}(x',\,x))}\; 
d\mu_{x}(v_{-})\,d\mu_{x}(v_{+})\,ds 
\\ &
=\Big(\frac{t'\,d_{\rm Cyg}(x,v_-)^{4}\,t'\,d_{\rm Cyg}(x,v_+)^{4}\,t^2}
{t\,d_{\rm Cyg}(x',v_-)^{4}\,t\,d_{\rm Cyg}(x',v_+)^{4}\;
d_{\rm Cyg}(x,v_-)^{4}\,d_{\rm Cyg}(x,v_+)^{4}}\Big)^{n}\;
d\lambda_{2n-1}(v_{-})\,d\lambda_{2n-1}(v_{+})\,ds
\\ &
=\frac{1}{d_{\rm Cyg}(v_-,v_+)^{4n}}\;
d\lambda_{2n-1}(v_{-})\,d\lambda_{2n-1}(v_{+})\,ds\,.
\end{align*}

\smallskip\noindent {\em (iii)~} Recall that the Liouville measure
$\vol_{T^1\hnc}$ (which is the Riemannian measure for Sasaki's metric
on $T^1\hnc$) disintegrates under the fibration $\pi:T^1\hnc \ra \hnc$
over the Riemannian measure $\vol_{\hnc}$ of $\hnc$, with conditional
measures the spherical measures on the unit tangent spheres:
$$
d\vol_{T^1\hnc}(v)=\int_{x\in \hnc} \dvol_{T^1_x\hnc}(v)\dvol_{\hnc}(x)\;.
$$ 
Let $x\in \hnc$ with last horospherical coordinate $t$. Since the
group $I_x$ of isometries of $\hnc$ fixing $x$ acts transitively on
$T^1_x\hnc$, and both $\mu_x$ and the Riemannian measure
$\vol_{T^1_x\hnc}$ are invariant under $I_x$, using the
$I_x$-equivariant homeomorphism $v\mapsto v_+$ from $T^1_x\hnc$ to
$\partial_\infty\hnc$, we have, for all $v\in T^1_x\hnc$ such that
$v_+\neq \infty$, using Assertion {\em (i)} for the last equality,
\begin{equation}\label{eq:relatvolfibunilambda}
d\vol_{T^1_x\hnc}(v)=\frac{\Vol(\SSS^{2n-1})}{\|\mu_x\|}\;d\mu_x(v_+)
=\frac{\Vol(\SSS^{2n-1})\;t^{n}}{\|\mu_x\|\,d_{\rm Cyg}(x,v_+)^{4n}}
\;d\lambda_{2n-1}(v_+)\,.
\end{equation}
By homogeneity and by Assertion {\em (i)}, we have, 
\begin{align*}
\|\mu_x\|&=\|\mu_{(0,0,1)}\|=
\int_{\CC^{n-1}\times \RR}\frac{d\xi\,dr}{((|\xi|^2+1)^2+r^2)^{n}}
\\ & = 2\Vol(\SSS^{2n-3})
\int_{[0,+\infty[^2}\frac{\rho^{2n-3}\;d\rho\,dr}{((\rho^2+1)^2+r^2)^{n}}\;.
\end{align*}
Using Mathematica for the first equation and the well known expression
of $\Ga(n+1/2)$ for the second one, we have
\begin{align*}
 2\int_{[0,+\infty[^2}\frac{s^{n-2}\,ds\,dr} {(s^2+r^2+2s+1)^{n}}
& =\frac{\sqrt{\pi}\;(n-2)!\;\Ga(n+1/2)}{(2n-2)!}\\
& =\frac{\sqrt{\pi}\;(n-2)!}{(2n-2)!}\;\frac{(2n)!\sqrt{\pi}}{2^{2n}\,n!}
=\frac{(2n-1)\,\pi}{2^{2n-1}\,(n-1)}\;.
\end{align*}
Using the change of variable $s=\rho^2$, since $\Vol(\SSS^{2n-1})=
\frac{\pi}{n-1}\,\Vol(\SSS^{2n-3})$, we hence have
$$
\|\mu_x\|=\frac{(2n-1)}{2^{2n-1}}\;\Vol(\SSS^{2n-1})\;.
$$
Hence, by Equations \eqref{eq:volhnc} and
\eqref{eq:relatvolfibunilambda}, using the homeomorphism $v\mapsto
(v_+,\pi(v)=(\zeta,u,t))$ from $T^1\hnc$ to $\partial_\infty \hnc
\times \hnc$, we have, for all $v\in T^1\hnc$ such that $v_+\neq
\infty$,
\begin{equation}\label{eq:liouville}
d\vol_{T^1\hnc}(v) =
\frac{2^{2n-3}}{(2n-1)\,t\,d_{\rm Cyg}((\zeta,u,t),v_+)^{4n}}
\;d\lambda_{2n-1}(v_+)\,d\zeta\,du\,dt\,.
\end{equation}

Consider the map $F:\Heis_{2n-1}\times\RR\,\ra\,\hnc$ defined by
$$
(\xi,r,s)\mapsto\Big(
\zeta=\frac{\xi}{1+(|\xi|^2-ir)e^{2s}},\;
u=-\Im\, \frac{|\xi|^2-ir}{1+(|\xi|^2-ir)e^{2s}},\;
t= \frac{e^{2s}(|\xi|^4+r^2)}{|1+(|\xi|^2-ir)e^{2s}|^2}\Big)\;.
$$
Note that $F(0,1,0)=(0,\frac{1}{2},\frac{1}{2})$. By the first two
centred formulas of \cite[page 113]{ParPau11MZ}, applied with the term
$t$ in these formulas equal to $2s+\ln 2$ and using Equation
\eqref{eq:horosphecoord}, noting that the sectional curvature in
\cite{ParPau11MZ} is normalised to be in $[-1,-\frac{1}{4}]$, the map
$s\mapsto F(\xi,r,s)$ is a geodesic line starting from $(\xi,r)$ and
ending in $(0,0)$. On this geodesic, $s$ and the time parameter in
Hopf's parametrisation differ only by an additive constant, hence have
the same differential.

Recall that by homogeneity, the two measures $\wt m_{\rm BM}$ and
$\vol_{T^1\hnc}$ are proportional. Hence, computing their (constant)
Radon-Nikodym derivative at $v\in T^1\hnc$ such that $v_-=(0,1)$ and
$\pi(v)=(0,\frac{1}{2},\frac{1}{2})$ (so that $v$ is tangent to the
geodesic line $s\mapsto F(0,1,s)$ at $s=0$, hence $v_+=(0,0)$), we
have, by Assertion {\em (ii)~} with $(\xi,r)$ parametrizing $v_-$ and
by Equation \eqref{eq:liouville},
\begin{align*}
\frac{d\vol_{T^1\hnc}}{d\wt m_{\rm BM}}&= 
\frac{2^{2n-3}\,d_{\rm Cyg}((0,1,0),(0,0,0))^{4n}} 
{(2n-1)\,\frac 12\,d_{\rm Cyg}((0,\frac{1}{2},\frac{1}{2}),(0,0,0))^{4n}}\;
\frac{d\zeta\,du\,dt}{d\xi\,dr\,ds}(0,1,0)\\ &
=\frac{2^{3n-2}}{(2n-1)}\;
\frac{d\zeta\,du\,dt}{d\xi\,dr\,ds}(0,1,0)
\,.
\end{align*}
The first claim of Assertion {\em (iii)} follows by computing the
Jacobian at $(0,1,0)$ of the map $F$, which is equal to
$\frac{1}{2^{n}}$, since at the point $(0,1,0)$, we have
$$
\frac{\partial \zeta}{\partial \xi}=
\frac{1}{1-i}\operatorname{Id}_{\CC^{n-1}},\;\;\;
\frac{\partial \zeta}{\partial r}=\frac{\partial \zeta}{\partial s}=0
,\;\;\;\frac{\partial u}{\partial r}=0,\;\;\;
\frac{\partial u}{\partial s}=-1,\;\;\;
\frac{\partial t}{\partial r}=\frac{1}{2}\,.
$$
The second claim follows from the facts that $\Vol(T^1M)=
\Vol(\SSS^{2n-1})\,\Vol (M)$ and that $\Vol(\SSS^{2n-1})=
\frac{2\,\pi^n}{(n-1)!}$.

\medskip\noindent {\em (iv)~} The first claim follows from Equation
\eqref{eq:skinninginfinity}.  By Equation \eqref{eq:volhorosphzetau},
we have
\begin{equation}\label{eq:volHun}
d\vol_{\partial \H_1}(\zeta,u,1)= \frac{1}{2}\;d\zeta\,du\,.
\end{equation}
Hence $\pi_*\wt \sigma_{\H_1}= 2\, \vol_{\partial \H_1}$, and by the
transitivity of the isometry group of $\hnc$ on the set of horoballs
in $\hnc$, the second claim of Assertion {\em (iv)} follows.
Therefore, by Lemma \ref{lem:horobvol},
$$
\|\sigma_{D^-}\|=\|\pi_*\sigma_{D^-}\|=
2\,\Vol(\Ga_{D^-}\bs\partial D^-)= 4n\,\Vol(\Ga_{D^-}\bs D^-)\,.
$$

\medskip\noindent {\em (v)~} By the transitivity of the isometry group
of $\hnc$ on the set of its geodesic lines, we may assume that $D^-$
is the geodesic line in $\hnc$ with points at infinity $(0,0)$ and
$\infty$. The map from the full-measure open subset
$\big\{(\zeta,u)\in\Heis_{2n-1}\,:\;\zeta\neq0\big\}$ in
$\Heis_{2n-1}$ to the product manifold $\SSS^{2n-3}\times\,]0,+\infty[
\,\times\, ]-\frac{\pi}{2}, \frac{\pi}{2}[$ defined by
\begin{equation}\label{eq:coordsigmarhotheta}
(\zeta,u)\mapsto \Big(\sigma=\frac{\zeta}{|\zeta|},\;
\rho=(|\zeta|^4+u^2)^{1/2},\;\theta=\arctan\frac{u}{|\zeta|^2}\Big)
\end{equation}
is a diffeomorphism. Since $|\zeta|=\sqrt{\rho\cos\theta}$ and
$u=\rho\sin\theta$, we have
\begin{align*}
d\zeta\,du=
\frac{1}{2} |\zeta|^{2n-4}\dvol_{\SSS^{2n-3}}\Big(\frac{\zeta}
{|\zeta|}\Big)\,d(|\zeta|^2)\,du
=\frac{1}{2} (\rho\cos\theta)^{n-2}\rho\dvol_{\SSS^{2n-3}}(\sigma)\,
d\rho\,d\theta\,.
\end{align*}
By Lemma \ref{lem:orthprojgeod}, by Equation \eqref{eq:rewritelem7} and
by Assertion {\em (i)}, using the homeomorphism sending
$v\in\normalout D^-$ to $v_+=(\zeta,u)\in \Heis_{2n-1} -\{(0,0)\}$ and
the inverse transformation $|\zeta|^2= \rho\cos\theta$,
$u=\rho\sin\theta$, we have
\begin{align*}
d\wt \sigma_{D^-}(v)&= 
e^{-2n\,\beta_{(\zeta,u)}(\,(0,\,0,\,(|\zeta|^4+u^2)^{1/2}),\;(0,\,0,\,1)\,)}\;
d\mu_{(0,\,0,\,1)}(\zeta,u)\\ & =\frac{(|\zeta|^4+u^2)^{n/2}\;d\zeta \,du}
{d_{\rm Cyg}((0,\,0,\,(|\zeta|^4+u^2)^{1/2}),(\zeta,u,0))^{4n}}=
\Big(\frac{(|\zeta|^4+u^2)^{1/2}}
{(|\zeta|^2+(|\zeta|^4+u^2)^{1/2})^2+u^2}\Big)^{n}\;d\zeta \,du
\\ &=
\frac{\cos^{n-2}\theta}{2^{n+1}\,(1+\cos\theta)^{n}}
\dvol_{\SSS^{2n-3}}(\sigma)\,\frac{d\rho}{\rho}\,d\theta\,.
\end{align*}
Thus,
\begin{equation}\label{eq:formskinningmeasgeod}
d\pi_*\wt \sigma_{D^-}(0,0,\rho)=
\frac{c'_n\,\Vol(\SSS^{2n-3})}{2^{n+1}}\;\frac{d\rho}{\rho}\;,
\end{equation}
where, using Mathematica, 
$$
c'_n=\int_{-\frac{\pi}{2}}^{\frac{\pi}{2}}\frac{\cos^{n-2}\theta}
{(1+\cos\theta)^{n}}\;d\theta
=\frac{2^{1-n}\;\sqrt{\pi}\;n!}{(n-1)\,\Ga(n+\frac{1}{2})}\,.
$$

The next step is to obtain a similar expression for the Riemannian
measure of the submanifold $\normalout D^-$ of $T^1\hnc$ (endowed with
Sasaki's metric). For every $x\in D^-$, let us denote by $\nu^1_x D^-$
the fiber over $x$ of the normal bundle map $v\mapsto \pi(v)$ from
$\normalout D^-$ to $D^-$. We endow $\nu^1_x D^-$ with the spherical
metric induced by the scalar product of the tangent space $T_x\hnc$ at
$x$. The Riemannian measure of $\normalout D^-$ disintegrates under
this fibration over the Riemannian measure of $D^-$ as
$$
d\vol_{\normalout D^-}(v)=\int_{x\in D^-}d\vol_{\nu^1_x D^-}(v)\dvol_{D^-}(x)\,.
$$
By looking at the expression \eqref{eq:riemethorosphecoord} of the
Riemannian metric of $\hnc$ in horospherical coordinates, using the
homeomorphism $\rho \mapsto x=(0,0,\rho)$ from $[0,+\infty[$ to $D^-$,
we have
$$
d\vol_{D^-}(x)=\frac{d\rho}{2\rho}\,.
$$
Hence
\begin{equation}\label{eq:formvolnormbundmeasgeod}
d\pi_*\vol_{\normalout D^-}(0,0,\rho)=\Vol(\SSS^{2n-2})\;\frac{d\rho}{2\rho}\,.
\end{equation}
We have $\Vol(\SSS^{2n-2})=\frac{2\;\pi^{n-\frac{1}{2}}}
{\Ga(n-\frac{1}{2})}=\frac{2^{2n-1}\,\pi^{n-1}\,(n-1)!}{(2n-2)!}$ and
$\Vol(\SSS^{2n-3})=\frac{2\;\pi^{n-1}} {(n-2)!}$.  Since
$\Ga(x+1)=x\Ga(x)$ for all $x>0$, Equations
\eqref{eq:formskinningmeasgeod} and \eqref{eq:formvolnormbundmeasgeod}
give the first claim of Assertion {\em (v)}. 

The second one follows, since pushforwards of measures preserves their
total mass, and since $\Vol(\Ga_{\normalout D^-}\bs \normalout D^-)=
\frac{\Vol(\SSS^{2n-2})}{m}\;\Vol(\Ga_{D^-}\bs D^-)$.

\bigskip\noindent {\em (vii)~} By the transitivity of the isometry
group of $\hnc$ on the set of its complex geodesic lines, we may
assume that $D^-$ is the complex geodesic line
$C=\{(w_0,w)\in\hnc\;:\;w=0\}$ or, in horospherical coordinates,
$C=\{(\zeta,u,t)\in \hnc\;:\; \zeta =0\}$.  By
\cite[p.~157]{Goldman99}, the orthogonal projection of $\hnc$ to $C$
is the map $(w_0,w)\mapsto (w_0,0)$, which in horospherical
coordinates on $\partial_\infty \hnc- \partial_\infty C$ extends as
$(\zeta,u,0)\mapsto (0,u,|\zeta|^2)$.

Hence, using the homeomorphism from $\normalout C$ to $\{(\zeta,u)
\in\Heis_{2n-1}\;:\;\zeta\neq 0\}$ sending a normal unit vector $v$ to
its point at infinity $v_+=(\zeta,u)$, we have, by Equation
\eqref{eq:rewritelem7} and by Assertion {\em (i)}, \ref{lem:busemann},
\begin{align*}
d\wt \sigma_{C}(v)&= 
e^{-2n\,\beta_{(\zeta,u)}(\,(0,\,u,\,|\zeta|^2),\;(0,\,0,\,1)\,)}\;
d\mu_{(0,\,0,\,1)}(\zeta,u)=\frac{1}{4^n|\zeta|^{2n}}\;d\zeta\,du\\ &=
\frac{1}{2^{2n-1}\,|\zeta|^{4}}\dvol_{\SSS^{2n-3}}
\big(\frac{\zeta}{|\zeta|}\big)\,du\,d(|\zeta|^2)\,.
\end{align*}
In particular,
$$
d\pi_*\wt \sigma_{C}(0,u,|\zeta|^2) =
\frac{\Vol(\SSS^{2n-3})}{2^{2n-1}}\;du\;\frac{d(|\zeta|^2)}{|\zeta|^{4}}\;.
$$

For every $x\in C$, let us denote by $\nu^1_x C$ the fiber over $x$ of
the normal bundle map $v\mapsto \pi(v)$ from $\normalout C$ to $C$,
endowed with the spherical metric induced by the scalar product of the
tangent space $T_x\hnc$ at $x$. The Riemannian measure of $\normalout
C$ disintegrates under this fibration over the Riemannian measure of
$C$ as
$$
d\vol_{\normalout C}(v)=\int_{x\in C}d\vol_{\nu^1_x C}(v)\dvol_{C}(x)\,.
$$
Using the homeomorphism $(u,t=|\zeta|^2) \mapsto x=(0,u,t)$ from
$\RR\times[0,+\infty[$ to $C$, and Equation
\eqref{eq:riemethorosphecoord}, we have
$$
d\vol_{C}(x)=\frac{du\,dt}{4\,t^2}=\frac{du\,d(|\zeta|^2)}{4\,|\zeta|^4}\,.
$$
Hence
$$
d\pi_*\vol_{\normalout C} (x)=\Vol(\SSS^{2n-3})\dvol_{C}(x)=
\frac{\Vol(\SSS^{2n-3})}{4}\;du\;\frac{d(|\zeta|^2)}{|\zeta|^{4}}\;.
$$
The result follows as in the end of the proof of the previous
Assertion. \cqfd

\medskip By Theorem \ref{theo:mainequicountdown}, we then have the
following counting and equidistribution result of common
perpendiculars. We first define some constants. 

Let $D^-$ be a horoball in $\hnc$ centred at a parabolic fixed point
of a lattice $\Ga$ in $\Isom(\hnc)$. If $D^+$ is also a horoball in
$\hnc$ centred at a parabolic fixed point of $\Ga$, let
$$
c(D^-,D^+)=\frac{4^{n}\,n!}{(2n-1)\,\pi^n}\;
\frac{\Vol(\Ga_{D^-}\bs D^-)\Vol(\Ga_{D^+}\bs D^+)}{\Vol(M)}\,.
$$
If $D^+$ is a geodesic line in $\hnc$ such that $\Ga_{D^+}\bs D^+$ is
compact, let
$$
c(D^-,D^+)=\frac{4^{n}\,(n!)^2}{(2n)!\,(2n-1)\,\pi}
\frac{\Vol(\Ga_{D^-}\bs D^-)\Vol(\Ga_{D^+}\bs D^+)}
{\Vol(M)}\,.
$$
If $D^+$ is a complex geodesic line in $\hnc$ such that $\Ga_{D^+}\bs
D^+$ has finite volume, let
$$
c(D^-,D^+)=\frac{4\,(n-1)}{(2n-1)\,\pi}
\frac{\Vol(\Ga_{D^-}\bs D^-)\Vol(\Ga_{D^+}\bs D^+)}
{\Vol(M)}\,.
$$

\bcoro\label{coro:complexhyperbo} Let $\Ga$ be a discrete group of
isometries of $\hnc$ such that the orbifold $M=\Ga\bs\hnc$ has finite
volume.  In each of the above three cases, if $m^+$ is the cardinality
of the pointwise stabiliser of $D^+$ in $\Ga$, then
$$
\N_{D^-,\,D^+}(t)\sim \frac{c(D^-,D^+)}{m^+}\;e^{2n\,t}\,.
$$
If $\Ga$ is arithmetic, then there exists $\kappa>0$ such that, as
$t\ra+\infty$,
$$
\N_{D^-,\,D^+}(t)=\frac{c(D^-,D^+)}{m^+}\;
e^{2n\,t}\big(1+\operatorname{O}(e^{-\kappa t})\big)\;.
$$
Furthermore, if $D^-$ is a horoball centred at a parabolic fixed
point, then the origins of the common perpendiculars from $D^-$ to the
images of $D^+$ under the elements of $\Ga$ equidistribute in
$\partial D^-$ to the induced Riemannian measure:\ 
\begin{equation}\label{eq:distribhorobcomplexhyp}
\frac{2\,n\,m^+\,\Vol(\Ga_{D^-}\bs D^-)}{c(D^-,D^+)}
\;e^{-2n\,t}\;\sum_{x\in\partial D^-} m_{t}(x) \;
\Delta_{x}\;\weakstar\; \vol_{\partial D^-}\,.
\;\;\;\Box
\end{equation}
\ecoro

For smooth functions $\psi$ with compact support on $\partial D^-$,
there is an error term in the equidistribution claim
\eqref{eq:distribhorobcomplexhyp} when the measures on both sides are
evaluated on $\psi$, of the form $\bigO(e^{-\kappa t}\,\|\psi\|_\ell)$
where $\kappa>0$ and $\|\psi\|_\ell$ is the Sobolev norm of $\psi$ for
some $\ell\in\NN$, by the remark following Theorem
\ref{theo:mainequicountdown}.

\section{A Mertens' formula for Heisenberg groups}
\label{sect:mertensheisenberg}

Let $K$ be an imaginary quadratic number field. We will denote, in
Sections \ref{sect:mertensheisenberg} to \ref{sect:chains}, by
$\OOO_K$ its ring of integers, by $D_K$ its discriminant, by
$\zeta_K$ its zeta function, and by $|\OOO_K^\times|$ the order of
the unit group of $\OOO_K$. Let $\tr,\n:K\ra \QQ$ be the (absolute)
trace and norm of $K$, that is $\tr (z)=z+\overline{z}=2\,\Re\;z$ and
$\n(z)=z\,\overline{z}=|z|^2$. We denote by $\langle
a,\alpha,c\rangle$ the ideal of $\OOO_K$ generated by $a, \alpha,
c\in\OOO_K$.

Let $\mmm$ be a nonzero ideal in $\OOO_K$.  We endow the ring
$\OOO_K/\mmm$ with the involution induced by the complex
conjugation. Let $\SU_q(\OOO_K/\mmm)$ be the finite group of
$3\times3$ matrices in $\OOO_K/\mmm$, having determinant $1$ and
preserving the Hermitian form $-z_0\overline{z_2}- z_2\overline{z_0} +
z_1\overline{z_1}$ on $(\OOO_K/\mmm)^3$. Let $B_q(\OOO_K/\mmm)$ be its
upper triangular subgroup.  The action by shears on $\OOO_K\times
\OOO_K \times \OOO_K$ of the nilpotent group $\Heis_3(\OOO_K)$ defined
in the Introduction preserves $\OOO_K\times\mmm\times\mmm$. In this
section, we will study the asymptotic of the counting function
$\Psi_{\mmm}$, where, for every $s\geq 0$, $\Psi_{\mmm}(s)$ is the
cardinality of
$$
\Heis_3(\OOO_K)\bs
\big\{(a,\alpha,c)\in\OOO_K\times\mmm\times\mmm\;:\;
\tr(a\,\overline{c})=\n(\alpha),\;\langle a,\alpha,c\rangle=\OOO_K,
\;\n(c)\leq s\big\}\,.
$$
When $\mmm=\OOO_K$, this map $\Psi_{\mmm}$ is the counting function,
in terms of their standard heights, of the rational points over $K$ in
the complex projective plane $\PP_2(\CC)$, that lie in Segre's
hyperconic with equation $2\,\Re\; u-|v|^2=0$ in the standard affine
chart with coordinates $(u,v)$.

\btheo\label{theo:countHeis} As $s\ra+\infty$, we have
$$
\Psi_{\mmm}(s)= \frac{\zeta(3)\,|B_q(\OOO_K/\mmm)|}
{2\,\pi\,|D_K|^{\frac{1}{2}}\,\zeta_K(3)\,|\SU_q(\OOO_K/\mmm)|}
\;s^2+\bigO(s^{2-\kappa})\,.
$$
\etheo

The particular case $\mmm=\OOO_K$ gives Theorem
\ref{theo:countHeisintro} in the Introduction. We will prove this
result simultaneously with the next one. We endow the $3$-dimensional
Heisenberg group $\Heis_{3}$ with its Haar  measure $\haarheis$ as in
the introduction. The following result is an equidistribution result
of the set of $\QQ$-points (satisfying some congruence properties) in
$\Heis_{3}$.  The particular case $\mmm=\OOO_K$ gives Theorem
\ref{theo:equidisHeisintro} in the introduction.

\btheo\label{theo:equidisHeis} As $s\ra+\infty$, we have
$$
\frac{\pi\,|D_K|^{\frac{3}{2}}\,\zeta_K(3)
\,|\SU_q(\OOO_K/\mmm)|}{\zeta(3)\,|B_q(\OOO_K/\mmm)|}\;s^{-2}
\sum_{\substack{(a,\,\alpha,\,c)\in \OOO_K\times\mmm\times\mmm,
\;0<\n(c)\leq s \\ \tr(a\,\overline{c})=\n(\alpha),
\;\langle a,\,\alpha,\,c\rangle=\OOO_K}}\;
\Delta_{(\frac{a}{c},\frac{\alpha}{c})}\;\weakstar\;\haarheis\,.
$$ 
\etheo

As said after Corollary \ref{coro:complexhyperbo}, for smooth
functions $\psi$ with compact support on $\Heis_3$, there is an error
term in this equidistribution result when the measures on
both sides are evaluated on $\psi$, of the form
$\bigO(s^{-\kappa}\,\|\psi\|_\ell)$ where $\kappa>0$ and
$\|\psi\|_\ell$ is the Sobolev norm of $\psi$ for some $\ell\in\NN$.

\medskip
\noindent{\bf Proof of Theorem \ref{theo:countHeis} and Theorem
  \ref{theo:equidisHeis}.}  As a preliminary remark, the
$3$-dimensional Heisen\-berg group $\Heis_{3}$ defined above contains
$\Heis_3(\OOO_K)$ (by the definition of the norm and trace of $K$), as a
(uniform) lattice, and identifies with the Heisenberg group
$\Heis_{3}$ defined in Section \ref{sec:cxhyp} by the change
of variable
$$
(\zeta,u)\mapsto (w_0=\frac{1}{2}(|\zeta|^2-iu),\;w=\zeta)\;,
$$
so that the Haar measures $\haarheis$ and $\lambda_3$ satisfy
\begin{equation}\label{eq:relatlambdatroihaarheis}
\lambda_3=2\;\haarheis\,.
\end{equation}

Let $q$ be the Hermitian form $-z_0\overline{z_2}- z_2\overline{z_0} +
z_1\overline{z_1}$ of signature $(2,1)$ on $\CC\times\CC\times\CC$
with coordinates $(z_0,z_1,z_2)$ (which is, up to isomorphism, the
unique Hermitian form over $K$ with signature $(2,1)$ and determinant
$-1$, see \cite[Ch.~10]{Scharlau85} for this cultural remark).  As
previously, we denote by $\SU_q$ the special unitary group of $q$. Let
$\Ga=\SU_q\cap\,\SL_3(\OOO_K)$ be the {\it Picard modular group} of
$K$, which is a nonuniform arithmetic lattice in $\SU_q$ by a result
of Borel and Harish-Chandra (see for instance \cite[\S
6.3]{ParPau10GT}). As another cultural remark, every nonuniform arithmetic
lattice in $\SU_q$ is commensurable to a Picard modular group (see for
instance \cite[\S~3.1]{Stover11}).

Consider the map from $\Heis_{3}$ to $\SU_q$ defined, in the two sets of
coordinates of $\Heis_{3}$ defined in the Introduction and in Section
\ref{sec:cxhyp}, by
$$
(w_0,w)\mapsto\left(\!\!\begin{array}{ccc}
1 & \overline{w} & w_0\\ 
0 & 1 & w\\ 0 & 0 & 1
\end{array}\!\!\right)\;\;\;{\rm or}\;\;\;
(\zeta,u)\mapsto\left(\!\!\begin{array}{ccc}
1 & \overline{\zeta} & \frac{|\zeta|^2}{2}-i\,\frac{u}{2}\\ 
0 & 1 & \zeta\\ 0 & 0 & 1
\end{array}\!\!\right)\;.
$$
This map is a Lie group isomorphism onto its image, by which we
identify from now on $\Heis_{3}$ and its image. Note that
$\Heis_{3}\cap\,\Ga$ is then exactly $\Heis_3(\OOO_K)$.

We denote by $\Ga_\mmm$ the Hecke congruence subgroup of
$\Ga$ modulo $\mmm$, that is the preimage, by the group
morphism $\Ga\ra \SL_3(\OOO_K/\mmm)$ of reduction modulo
$\mmm$, of the upper triangular subgroup of $\SL_3(\OOO_K/\mmm)$.

As previously, we denote by $\PSU_q$ the quotient Lie group
$\SU_q/\{\id,j\id,j^2\id\}$ where $j=e^{\frac{2i\pi}{3}}$. For every
subgroup $G$ of $\SU_q$, we denote by $\overline{G}$ its image in
$\PSU_q$, and again by $g$ the image in $\PSU_q$ of any element $g$ of
$\SU_q$.

We denote by $\Ga_{\H_1}$ the stabiliser in $\Ga_\mmm$ of the horoball
$$
\H_1=\{[w_0:w:1]\in \PP_2(\CC)\;:\; 2\,\Re\,w_0-|w|^2\geq 1\}
$$
centred at $\infty= [1:0:0]$ in the projective model of $\hdc$. Note
that this agrees with our notation for horoballs introduced in Section
\ref{sec:cxhyp}.  The group $\Ga_{\H_1}$ is independent of $\mmm$,
since $\Ga_{\mmm}$ is Hecke's congruence subgroup of $\Ga$.  Note that
$\Heis_3(\OOO_K)$ is contained in $\Ga_{\H_1}$, and that the
projection map from $\Heis_3(\OOO_K)$ to $\overline{\Heis_3(\OOO_K)}$
is injective.

Let $B_q$ be as defined in Section \ref{sec:cxhyp}. Note that
$\Ga_{\H_1}= B_q\cap\Ga_\mmm= B_q\cap\Ga$, since an isometry
of $\hdc$ fixes $\infty$ if and only if it preserves $\H_1$.  We claim
that the index of $\Heis_3(\OOO_K)$ in $\Ga_{\H_1}$ is
\begin{equation}\label{eq:calcindexheisGainfty}
[\Ga_{\H_1}:\Heis_3(\OOO_K)]=|\OOO_K^\times|\,.
\end{equation}
Indeed, it is the cardinality of the set of $(a_1,a_2,a_3)\in
(\OOO_K^\times)^3$ such that $a_1a_2a_3=1, \,a_3\,\overline{a_1}=1,
\;|a_2|=1$, which is the cardinality of the set of $(a_1,a_2) \in
(\OOO_K^\times)^2$ such that ${a_1}^2=1/a_2$. A separate treatment of
the cases $D_K=-3,-4$ and of the general case gives the result.  A
similar argument shows that the projection map from $\Ga_{\H_1}$ to
$\overline{\Ga_{\H_1}}$ is $(1+2\,\delta_{D_K,-3})$-to-$1$ (in particular
injective when $D_K\neq-3$).

For all $a,\alpha,c\in\OOO_K$ and $\lambda\in\CC$, we have $\langle
\lambda a, \lambda \alpha,\lambda c\rangle= \langle a,\alpha, c
\rangle =\OOO_K$ if and only if $\lambda\in \OOO_K^\times$. Therefore
the cardinality of the fibers of the projection map from
$\{(a,\alpha,c)\in \OOO_K\times\mmm \times \mmm\;:\; \langle
a,\alpha,c\rangle=\OOO_K\}$ to $\PP_2(\CC)$ is equal to
$|\OOO_K^\times|$, and no two distinct elements of such a fiber are
sent one to the other by an element of $\Heis_3(\OOO_K)$. By
\cite[Prop.~6.5 (2)]{ParPau10GT}, the orbit
$\overline{\Ga_\mmm}\cdot\infty$ of $\infty=[1:0:0]$ under
$\overline{\Ga_\mmm}$ is equal to the set of elements of $\PP_2(\CC)$
which can be written in homogeneous coordinates $[a:\alpha:c]$ where
$(a,\alpha,c)\in \OOO_K\times\mmm\times\mmm$, $\langle
a,\alpha,c\rangle=\OOO_K$ and $\tr(a\,\overline{c})= \n(\alpha)$. 

Let $g\in \SU_q$ be such that $g\H_1$ and $\H_1$ are disjoint (there
are only finitely many double classes $[g]\in \overline{\Ga_{\H_1}}
\,\bs \overline{\Ga_\mmm}/ \,\overline{\Ga_{\H_1}}$ for which this is
not the case). If $\left(\!\!\begin{array}{c} a_g \\ \alpha_g \\
    c_g \end{array} \!\!\right)$ is the first column of $g$, then
$g\cdot \infty = [a_g:\alpha_g:c_g]$. Futhermore, by
\cite[Lem.~6.3]{ParPau10GT} and since the sectional curvature is
normalised to have maximum $-1$, the length of the common
perpendicular $ \delta_g$ between $g\,\H_1$ and $\H_1$ is then
\begin{equation}\label{eq:computlongcomperp}
\ell(\delta_g)=\ln |c_g|-\ln 2\;.
\end{equation}

We use, in the last one of the following equalities, Equation
\eqref{eq:calcindexheisGainfty} and Corollary
\ref{coro:complexhyperbo} with $n=2$, $\Ga=\overline{\Ga_\mmm}$ and
$D^-=D^+=\H_1$ (whose pointwise stabilisers in $\Ga_\mmm$ are
trivial). We hence have, for some $\kappa>0$,
\begin{align}
\Psi_\mmm(s)& = |\OOO_K^\times|\;\;\card\;\;_{\overline{\Heis_3(\OOO_K)}}\,\bs
\Big\{[a:\alpha:c]\in\PP_2(\CC)\;:\;
\begin{array}{c}(a,\alpha,c)\in\OOO_K\times\mmm\times\mmm,\\
\langle a,\alpha,c\rangle=\OOO_K,\\
\tr(a\,\overline{c})=\n(\alpha),\;\n(c)\leq s\end{array}
\Big\}\nonumber\\ 
&=
|\OOO_K^\times|\;\;\card\;\;_{\overline{\Heis_3(\OOO_K)}}\,\bs
\big\{[a:\alpha:c]\in\overline{\Ga_\mmm}\cdot\infty\;:\;
\begin{array}{c}(a,\alpha,c)\in\OOO_K\times\mmm\times\mmm,\\
  \langle a,\alpha,c\rangle=\OOO_K,\;\n(c)\leq s\end{array}
\big\}\nonumber\\ 
&=
|\OOO_K^\times|\;\;\card\;
\big\{[g]\in\overline{\Heis_3(\OOO_K)}\,\bs
\overline{\Ga_\mmm}/\,\overline{\Ga_{\H_1}}\;:\;\n(c_g)\leq s
\big\}\nonumber\\ 
&=
|\OOO_K^\times|\;
[\,\overline{\Ga_{\H_1}}:\overline{\Heis_3(\OOO_K)}\,]\;\card
\big\{[g]\in\overline{\Ga_{\H_1}}\,\bs
\overline{\Ga_\mmm}/\,\overline{\Ga_{\H_1}}\;:\; \ell(\delta_g)\leq
\frac{\ln s}{2}-\ln 2 \big\}+\bigO(1)\nonumber\\ 
&=
|\OOO_K^\times|\;
\frac{[\Ga_{\H_1}:\Heis_3(\OOO_K)]}{1+2\,\delta_{D_K,-3}}
\;\;\N_{\H_1,\,\H_1}\big(\frac{\ln
  s}{2}-\ln 2\big)+\bigO(1) \nonumber \\ 
  &=
  \frac{2\;|\OOO_K^\times|^2\;\Vol(\,\overline{\Ga_{\H_1}}\,\bs
  \H_1)^2}
{3\,\pi^2\,(1+2\,\delta_{D_K,-3})\,
\Vol(\,\overline{\Ga_\mmm}\,\bs
  \hdc)} \;s^2\,(1+\bigO(s^{-\kappa}))\;.\label{eq:longcalcul}
\end{align}

\medskip The next two lemmas are devoted to the computation of the two
volumes that appear in the previous line. 

\blemm\label{lem:volumecuspHeis} We have 
$\displaystyle \Vol(\,\overline{\Ga_{\H_1}}\,\bs \H_1)=
\frac{(1+2\,\delta_{D_K,-3})\,|D_K|}{8\,|\OOO_K^\times|}\,.$ 
\elemm

\dem We follow arguments similar to the ones in the reference \cite[\S
4]{KimKim09}, which uses the same convention as \cite{Parker98} for
the Riemannian measure on $\Heis_{2n-1}$, see also \cite[\S
3.1]{FalPar06} when $D_K=-3$.  Let $n=2$. As in \cite{Parker98}, we
endow $\Heis_{3}$ with the left-invariant Riemannian metric
$$
(du+2\,\Im(d\zeta\;\overline{\zeta})\,)^2+
4\,d\zeta\;\overline{d\zeta}\,,
$$
whose Riemannian volume is $\vol_{\Heis_{3}}=4\lambda_3$. In
particular, the $\Heis_{3}$-equivariant map from $\partial\H_1$ to
$\Heis_{3}$, which in horospherical coordinates maps $(\zeta,u,t)$ to
$(\zeta,u)$, sends $\vol_{\partial\H_1}$ to $\frac{1}{2}\,\lambda_3=
\frac{1}{8}\,\vol_{\Heis_{3}}$, by Equation
\eqref{eq:riemethorosphecoord}.

Let $t_K$ be the minimal vertical translation in $\Heis_3(\OOO_K)$,
that is, the minimal $s>0$ such that $(w_0=\frac{is}{2},w=0) \in
\Heis_3(\OOO_K)$. In particular,
$$
t_K=\min\{s>0\;:\;\frac{is}{2}\in\OOO_K\}\;.
$$ 
Recall that $\OOO_K=\ZZ+\ZZ \frac{i\sqrt{|D_K|}}{2}$ if $D_k\equiv
0\mod 4$ and $\OOO_K=\ZZ+\ZZ\frac{1+ i\sqrt{|D_K|}}{2}$ if $D_k\equiv
1\mod 4$. Hence $t_K= \sqrt{|D_K|}$ if $D_k\equiv 0\mod 4$ and $t_K=
2\,\sqrt{|D_K|}$ if $D_k\equiv 1\mod 4$.

Consider the following set
\begin{align*}
\Pi(\Heis_3(\OOO_K))&=
\{w\in\CC\;:\;\exists\;w_0\in\CC,\;(w_0,w)\in \Heis_3(\OOO_K)\}
\\ &=\{w\in\OOO_K\;:\;\exists\;w_0\in\OOO_K,\;\tr (w_0)=\n(w)\}
=\n^{-1}(\tr(\OOO_K))\,.
\end{align*}
We have $\tr(\OOO_K)=\ZZ$ if $D_K\equiv 1 \mod 4$. If $D_K\equiv 0
\mod 4$, then $\tr(\OOO_K)=2\ZZ$ and for all $a,b\in\ZZ$, the integer
$\n\big(a+bi\,\frac{\sqrt{|D_K|}}2\big)=a^2+b^2\frac{|D_K|}4$ is even
if and only if
\begin{itemize}
\item  $a$ is even, when $\frac{|D_K|}4$ is even,  
\item $a-b$ is even, when $\frac{|D_K|}4$ is odd.
\end{itemize}
Hence $\Pi(\Heis_3(\OOO_K))$ is a $\ZZ$-sublattice in
$\OOO_K$ with index equal to $1$ if $D_K\equiv 1 \mod 4$ and equal to
$2$ if $D_K\equiv 0 \mod 4$.

By \cite[Lem.~1.2]{Parker98}, we have
$$
\Vol_{\Heis_{3}}(\Heis_3(\OOO_K)\bs\Heis_{3})= 
4\,t_K\,\Vol(\CC/\Pi(\Heis_3(\OOO_K)))\;.
$$
Note that $\Vol(\CC/\OOO_K)= \frac{\sqrt{|D_K|}}{2}$ since $\OOO_K$ is
generated as a $\ZZ$-module by $1$ and $(D_K+i\sqrt{|D_K|})/2$.  Hence
\begin{equation}\label{eq:parkercuspnilp}
\Vol_{\Heis_{3}}(\Heis_3(\OOO_K)\bs\Heis_{3})= 
4\,t_K\,[\OOO_K:\Pi(\Heis_3(\OOO_K))]\;\Vol(\CC/\OOO_K))=4\,|D_K|\;.
\end{equation}

As already seen, the map $\Ga_{\H_1}\ra \overline{\Ga_{\H_1 }}$ is a
$3$-to-$1$ map if $D_K=-3$ and injective otherwise.  By Lemma
\ref{lem:horobvol}, we then have
\begin{align*}
\Vol(\,\overline{\Ga_{\H_1}}\,\bs \H_1)
&=\frac{1}{4}\;\Vol(\,\overline{\Ga_{\H_1}}\,\bs \partial \H_1)\\ 
&=\frac{1}{4\,[\,\overline{\Ga_{\H_1}}:\overline{\Heis_3(\OOO_K)}\,]}\;
\Vol(\,\overline{\Heis_3(\OOO_K)}\,\bs \partial \H_1)\\ 
&=\frac{1+2\,\delta_{D_K,-3}}
{32\,[\Ga_{\H_1}:\Heis_3(\OOO_K)]}\;
\Vol_{\Heis_{3}}(\Heis_3(\OOO_K)\bs \Heis_3)\,.
\end{align*}
The result hence follows from Equation \eqref{eq:calcindexheisGainfty}
and Equation \eqref{eq:parkercuspnilp}. \cqfd

\medskip The volume of the orbifold $\overline{\Ga_\mmm}\,\bs \hdc$ is
known (see for instance \cite{Holzapfel89,Stover11}).  We only give a
proof since the normalisation of the measures is a bit tricky to
follow amongst the various references.

\blemm \label{lem:calcvolPicmod} We have
$\displaystyle \Vol(\,\overline{\Ga_\mmm}\,\bs \hdc)=
\frac{(1+2\,\delta_{D_K,-3})\,|D_K|^{5/2}\,\zeta_K(3)\,|\SU_q(\OOO_K/\mmm)|}
{48\,\pi\,\zeta(3)\,|B_q(\OOO_K/\mmm)|}\,.$ 
\elemm

\dem
In the complex ball $$\BB^2_\CC=\{(z_1=x_1+iy_1,z_2=x_2+iy_2)\in
\CC^2\;:\; |z_1|^2+|z_2|^2 < 1\}\,,$$ the Riemannian metric invariant
under its group of biholomorphisms, which makes it isometric to the
Siegel domain $\hdc$, has volume form at the point $(0,0)$ equal
to
$$
4\;dx_1\wedge dy_1\wedge dx_2\wedge dy_2
$$
(see \cite[page 105]{Goldman99}, which normalises the sectional
curvature to be in $[-1,-\frac{1}{4}]$, as we normalise it to be in
$[-4,-1]$ in this paper). The volume form $\vol_{\rm Hol}$ at the
point $(0,0)$ of $\BB^2_\CC$ used in \cite[page 86]{Holzapfel89} for
its volume computation is
$$
\frac{6}{\pi^2}\;dx_1\wedge dy_1\wedge dx_2\wedge dy_2\;.
$$
The Cayley transform from the Siegel domain to the ball model of the
complex hyperbolic space conjugates $\overline{\Ga}$ to the
Picard modular group $\Ga$ used by Holzapfel in loc. cit..  Hence
$\Vol(\,\overline{\Ga}\,\bs\hdc)= \frac{2\pi^2}{3}\;
\Vol_{\rm Hol} (\Ga\bs\BB^2_\CC)$. By the Main Theorem 4.9 of
\cite[page 83]{Holzapfel89}, we have
$$
\Vol_{\rm Hol}(\Ga\bs\BB^2_\CC)
=\frac{(1+2\,\delta_{D_K,-3})\,|D_K|^{5/2}\,L_K(3)}{32\,\pi^3}\;,
$$
where $L_K(s)=\sum_{n=1}^\infty \chi_K(n)/n^s$ is Dirichlet's
$L$-series of $K$, whose character is given by Jacobi's symbol
$\chi_K(n)= \big(\frac{D_K}{n}\big)$. Hence, using the well known
relation $L_K(s)=\zeta_K(s)/\zeta(s)$ between Dirichlet's
$L$-series and Dedekind's zeta function of $K$, the result follows
since 
$$
[\,\overline{\Ga}\,:\,\overline{\Ga_\mmm}\,]
=[\Ga:\Ga_\mmm]=[\SU_q(\OOO_K/\mmm):B_q(\OOO_K/\mmm)]=
\frac{|\SU_q(\OOO_K/\mmm)|}{|B_q(\OOO_K/\mmm)|}\,. 
$$
\cqfd

\medskip Theorem \ref{theo:countHeis} follows from Lemma
\ref{lem:calcvolPicmod}, Lemma \ref{lem:volumecuspHeis} and Equation
\eqref{eq:longcalcul}.  

\medskip Let us prove now Theorem \ref{theo:equidisHeis}. 

The orthogonal projection map $f:\partial_\infty\hdc-\{\infty\}\ra 
\partial \H_1$ is the homeomorphism defined by $[w_0:w:1]\mapsto
(\zeta=w, u=-2\,\Im\,w_0,1)$ using the homogeneous coordinates on
$\partial_\infty\hdc-\{\infty\}$ and the horospherical coordinates on
$\partial \H_1$ (see Equation \eqref{eq:horosphecoord}). Let
$x\in\partial \H_1$ be the origin of a common perpendicular of length
at most $t$ from $\H_1$ to an element $\ga\H_1$ for some
$\ga\in\Ga_\mmm$ not fixing $\infty$. By the previous arguments, $x$
is the orthogonal projection on $\H_1$ of the point at infinity of
this horoball $\ga\H_1$. This point at infinity may be written
$[\frac{a}{c}: \frac{\alpha}{c} :1]$ for some triple $(a,\,\alpha,\,c)
\in \OOO_K\times \mmm \times\mmm$ with $\langle a,\,\alpha,\,c\rangle
=\OOO_K$, $\tr(a\,\overline{c})=\n(\alpha)$ and $0<\n(c)\leq
4\,e^{2t}$ (using Equation \eqref{eq:computlongcomperp}). There are
exactly $|\OOO^\times_K|$ such triples. Hence by Equation
\eqref{eq:distribhorobcomplexhyp}, considering the value of
$C(D^-,D^+)$ for $D^-=D^+=\H_1$, using the horospherical coordinates
on $\partial\H_1$, we have, as $t\ra+\infty$,
\begin{equation}\label{eq:distribhorobPic}
\frac{3\,\pi^2\,\Vol(\,\overline{\Ga_\mmm}\,\bs \hdc)}
{8\,\Vol(\,\overline{\Ga_{\H_1}}\,\bs \H_1)\,|\OOO_K^\times|}
\;e^{-4\,t}\;\sum_{\substack{(a,\,\alpha,\,c)\in
\OOO_K\times\mmm\times\mmm,\;0<\n(c)\leq 4\,e^{2\,t} \\ 
\tr(a\,\overline{c}) = \n(\alpha), 
\;\langle a,\,\alpha,\,c\rangle=\OOO_K}}\;
\Delta_{(\frac{\alpha}{c},\,-2\,\Im\,\frac{a}{c},\,1)}\;
\weakstar\; \vol_{\partial \H_1}\,.
\end{equation}

The image of the Haar measure $\haarheis$ (defined in Equation
\eqref{eq:defhaarheis}) by $f$ is, by Equation
\eqref{eq:riemethorosphecoord},
$$
f_*\haarheis=\vol_{\partial \H_1}\;.
$$
Using the change of variables $s=4\,e^{2t}$, the identification of
$\partial_\infty\hdc-\{\infty\}$ with $\Heis_3$ and the continuity of
the pushforward by $f^{-1}$ of the measures on $\partial \H_1$ applied
to Equation \eqref{eq:distribhorobPic}, we hence have, as $s\ra
+\infty$,
$$
\frac{6\,\pi^2\,\Vol(\,\overline{\Ga_\mmm}\,\bs \hdc)}
{\Vol(\,\overline{\Ga_{\H_1}}\,\bs \H_1)\,|\OOO_K^\times|}
\;s^{-2}\;\sum_{\substack{(a,\,\alpha,\,c)\in
\OOO_K\times\mmm\times\mmm,
\;0<\n(c)\leq s \\ \tr(a\,\overline{c})=\n(\alpha),
\;\langle a,\,\alpha,\,c\rangle=\OOO_K}}\;
\Delta_{(\frac{a}{c},\frac{\alpha}{c})}\;\weakstar\;\haarheis\;.
$$
Finally, Theorem \ref{theo:equidisHeis} follows from this and from
Lemma \ref{lem:calcvolPicmod} and Lemma \ref{lem:volumecuspHeis}.
\cqfd

\medskip \rem A result of Feustel \cite{Feustel77} (see also
\cite[page 280]{Holzapfel98} and \cite{Zink79}) says that the map,
which associates to a parabolic fixed point of $\overline{\Ga}$ the
fractional ideal generated by its homogeneous coordinates in $\OOO_K$,
induces a bijection from the set of cusps (that is, of orbits under
$\overline{\Ga}$ of its parabolic fixed points) to the set of ideal
classes of $K$. In Theorem \ref{theo:countHeisintro} and Theorem
\ref{theo:equidisHeisintro}, replacing in its proof $D^+=\H_1$ by a
horoball centred at a parabolic fixed point $p$ of $\overline{\Ga}$
not in the orbit of $\infty$ (which hence exists if and only if $D_K
\neq -3, -4,-7,-8,-11,-19, -43,-67,-163$), we can obtain a counting
and equidistribution result with error term in $\Heis_{3}$ of the
points in $\Ga\cdot p$. But the volume of the quotient of this new
$D^+$ by its stabiliser in $\overline{\Ga}$ is not explicit for the
moment, hence we would not have results as precise as in the case $p=
\infty$.

\medskip \brema\label{rem:highdim} {\rm Theorem \ref{theo:countHeis}
  and Theorem \ref{theo:equidisHeis} have generalisations in higher
  dimension. Let $n\geq 2$, and let $(w,w')\mapsto w\cdot
  \overline{w'}$ be the standard Hermitian scalar product on
  $\CC^{n-1}$. Let 
$$ \Heis_{2n-1}=\{(w_0,w)\in\CC\times\CC^{n-1}
  \;:\; 2\,\Re\;w_0 = w\cdot \ov{w}\}\,, 
$$ 
with law $(w_0,w)(w'_0,w')= (w_0+w'_0+ w'\cdot\overline{w}, w+w')$, be
the Heisenberg group of dimension $2n-1$, which identifies with the
boundary at infinity of the Siegel domain $\hnc$ with $\infty$
removed. Let $q$ be the Hermitian form defined in Section
\ref{sec:cxhyp}. Let $\SU_q$ be its special unitary group and
$\Ga=\SU_q\cap\,\M_n(\OOO_K)$, which is an arithmetic lattice in
$\SU_q$. Then Corollary \ref{coro:complexhyperbo} (which is valid in
any dimension), applied with the image $\overline{\Ga}$ of $\Ga$ in
$\PSU_q$ and with $D^-=D^+$ the horoball of points in $\hnc$ with last
horospherical coordinates at least $1$, gives a counting and
equidistribution result with error term in $\Heis_{2n-1}$ of the
points in $\overline{\Ga}\cdot \infty- \{\infty\}$. The volume of
$\overline{\Ga}\bs \hnc$ could be computed using \cite{EmeSto13}, up
to computing the index of $\Ga$ in a principal arithmetic subgroup
containing it. But the volume of the cusp corresponding to $\infty$ in
$\overline{\Ga}\bs \hnc$ is not explicit for the moment, hence we
would not have results as precise as in the case $n=2$.

Other counting and equidistribution results of arithmetically defined
points in the Heisenberg group $\Heis_{2n-1}$ may be obtained by
varying the integral Hermitian form $q$ of signature $(1,n)$ and the
arithmetic lattice $\Ga$ in $\SU_q$.  }\erema

\section{Counting cubic points over quadratic imaginary 
fields in the projective plane}
\label{sect:cubicpoints}

Let $K$ be an imaginary quadratic number field.  Let $q$ be the
Hermitian form $-z_0\overline{z_2} -z_2\overline{z_0}+|z_1|^2$ on
$\CC^3$ (the following result could be adapted to any Hermitian form
on $K^3$ with complex signature $(1,2)$, see for instance
\cite[Ex.~1.6(iv), p.~351]{Scharlau85} for their classification).  We
will say that a point in the complex projective plane $\PP_2(\CC)$ is
{\it isotropic} (respectively that two projective points are {\it
  orthogonal}) if the corresponding complex lines in $\CC^3$ are
isotropic (respectively orthogonal) for $q$.

The Galois group $\Gal(\CC|K)$ acts naturally on $\PP_2(\CC)$ by
$\sigma[z_0:z_1:z_2]= [\sigma z_0:\sigma z_1:\sigma z_2]$ using
homogeneous coordinates. A point $z\in \PP_2(\CC)$ will be called {\it
  Hermitian cubic} over $K$ if it is cubic over $K$ (that is, if its
orbit under $\Gal(\CC|K)$ has exactly three points), and if its other
conjugates $z',z''$ over $K$ are isotropic and orthogonal to $z$.

We will denote by $(\ga, z)\mapsto \ga\cdot z$ the projective action
of $\SL_3(\CC)$ or $\PSL_3(\CC)$ on $\PP_2(\CC)$. Recall that $\SU_q$
is the real Lie group of linear automorphisms of $\CC^3$ having
determinant $1$ and preserving $q$. Let $\Ga=\SU_q(\OOO_K)=\SU_q \cap
\SL_3(\OOO_K)$, which is an arithmetic lattice in $\SU_q$. The action
of $\Ga$ on $\PP_2(\CC)$ preserves the sets of rational, cubic and
Hermitian cubic points of $\PP_2(\CC)$ over $K$. The number of
orbits of Hermitian cubic points is infinite, which explains why we
will restrict ourselves to a given orbit below. Let $\Ga_\infty$ be
the stabiliser of $\infty=[1:0:0]$ in $\Ga$, which preserves the Cygan
distance $d_{\rm Cyg}$ and its modifications $d'_{\rm Cyg},d''_{\rm
  Cyg}$ on $\Heis_{3}=\{[w_0:w:1]\in\PP_2(\CC)\;:\; 2\,\Re\; w_0 =
|w|^2\}$.

Let $z\in \PP_2(\CC)$ be Hermitian cubic over $K$, and let $z',z''$ be
its other conjugates over $K$. Since $z',z''$ are distinct, irrational
over $K$, and isotropic for $q$, they lie in $\Heis_{3}$, and we
may define the {\it complexity} $c(z)$ of $z$ as the inverse of a
modification of the Cygan distance between its conjugates over $K$:
$$
c(z)=d''_{\rm Cyg}(z',z'')^{-1}\,.
$$
The complexity of $z$ is in particular invariant under $\Ga_\infty$.
It is part of the proof of the following result that the number of
projective points, that are Hermitian cubic over $K$, belong to a
given orbit of $\Ga=\SU_q(\OOO_K)$, and have complexity at most $s$,
is finite for every $s\geq 0$. It turns out that the use of $d''_{\rm
  Cyg}$ instead of the actual Cygan distance $d_{\rm Cyg}$ allows in
this Section to give precise asymptotic results with error terms,
instead of upper/lower estimates that differ by a multiplicative
constant.

Let us now state and prove an asymptotic estimate as $s\ra+\infty$ of
the counting function of the set of points in the projective plane
that are Hermitian cubic over $K$, in a given orbit of for instance a
congruence subgroup of $\SU_q(\OOO_K)$, with complexity at most $s$.

\btheo\label{theo:countingcubic} Let $K$ be an imaginary quadratic
number field. Let $z_0\in\PP_2(\CC)$ be Hermitian cubic over $K$. Let
$G$ be a finite index subgroup of $\PSU_q(\OOO_K)$, and let $G_
\infty$ be the stabiliser of $\infty=[1:0:0]$ in $G$. Then there
exists $\kappa>0$ such that, as $s\ra+\infty$,
\begin{multline*}
\card\{z\in G_\infty\bs G\cdot z_0\;:\; c(z)\leq s\}\\=
\frac{256\,\ln|\lambda_0|\;\zeta(3)\;[\PSU_q(\OOO_K)_\infty:G_\infty]}
{3\,\iota_0\,n_0\,|\OOO_K^\times|\,|D_K|^{3/2}\,\zeta_K(3)\,
[\PSU_q(\OOO_K):G]}\;s^4\,(1+\bigO(s^{-\kappa}))\;,
\end{multline*}
where $|\lambda_0|$ is the smallest modulus $>1$ of an
eigenvalue of an element in $G$ fixing the other two conjugates
$z'_0$, $z''_0$ of $z_0$ over $K$, where $\iota_0=2$ if $G$ contains
an element exchanging $z'_0$ and $z''_0$ and $\iota_0=1$ otherwise,
and where $n_0$ is the cardinality of the pointwise stabiliser in $G$
of the projective line through $z'_0$ and $z''_0$.  
\etheo

\dem The group $G$, which has finite index in the arithmetic lattice
$\overline{\Ga}=\PSU_q(\OOO_K)$, is also an arithmetic discrete group
of isometries with finite covolume in the complex hyperbolic plane
$\hdc$. By Lemma \ref{lem:calcvolPicmod}, we have
\begin{align}
\Vol(G\bs \hdc)&=
[\,\overline{\Ga}\,:G]\;\Vol(\,\overline{\Ga}\,\bs \hdc)
\nonumber \\ & =
\frac{[\,\overline{\Ga}\,:G]\,(1+2\,\delta_{D_K,-3})\,
|D_K|^{5/2}\,\zeta_K(3)}{48\,\pi\,\zeta(3)}\,.\label{eq:covolG}
\end{align}
Let $n_\infty=[\,\overline{\Ga}_\infty:G_\infty]$ be the index of
$G_\infty$ in the stabiliser of $\infty$ in $\overline{\Ga}$.
By Lemma \ref{lem:volumecuspHeis}, since $G_\infty$ is equal to the
stabiliser $G_{\H_1}$ in $G$ of the horosphere $\H_1$, we have
\begin{equation}\label{eq:covolGinfty}
\Vol(G_\infty\bs \H_1)=\Vol(G_{\H_1}\bs \H_1)=
\frac{n_\infty\,(1+2\,\delta_{D_K,-3})\,|D_K|}
{2\,|\OOO_K^\times|}\,.
\end{equation}

\bprop\label{prop:caracfixpoint} A point $z_0\in \PP_2(\CC)$ is
Hermitian cubic over $K$ if and only if there exists $\ga_0\in
\PSU_q(\OOO_K)$ of infinite order and $K$-irreducible such that $z_0$
is the only fixed point of $\ga_0$ that belongs to the positive cone
of $q$ in $\PP_2(\CC)$. 
\eprop

Recall that an element of $\PSU_q(K)$ is {\it $K$-irreducible} if it
does not preserve a point or a line defined over $K$ in $\PP_2(\CC)$.

\medskip 
\dem (Y.~Benoist) Assume first that $\ga_0\in \PSU_q(\OOO_K)$
has infinite order and is $K$-irreducible. Then $\ga_0$ is not an
elliptic element, since elliptic elements of $\PSU_q(\OOO_K)$ have
finite order. It is not parabolic, since the fixed points in
$\partial_\infty\hdc$ of the parabolic elements of $\PSU_q(\OOO_K)$
are rational over $K$ (see \cite{Holzapfel78} or \cite[page
290]{Holzapfel98}).  Hence, it is loxodromic, and fixes exactly two
distinct points $z'_0$ and $z''_0$ in $\partial_ \infty\hdc$. In
particular, $z'_0$ and $z''_0$ are isotropic for $q$. Since $\ga_0$
belongs to $\PSL_3(\CC)$ and preserves $\partial_ \infty\hdc$, it
preserves the unique complex projective lines $L'$ and $L''$ tangent
to $\partial_ \infty\hdc$ at the points $z'_0$ and $z''_0$, 
respectively. Note that $L'$ and $L''$ are exactly the sets of points
in $\PP_2(\CC)$ which are orthogonal to $z'_0$ and $z''_0$, 
respectively. The projective lines $L'$ and $L''$ meet at exactly one
point $z_0$, which belongs to the complement in $\PP_2(\CC)$ of
$\hdc\cup\partial_ \infty\hdc$. This complement is exactly the
positive cone of $q$ in $\PP_2(\CC)$. The fixed points $z_0, z'_0,
z''_0$ of $\ga_0$ are at most cubic over $K$, since $\ga_0$ has
coefficients in $K$. They are exactly cubic and conjugates, since
$\ga_0$ is $K$-irreducible. Hence $z_0$ is Hermitian cubic, and is the
only fixed point of $\ga_0$ in the positive cone of $q$.

Conversely, let $z_0\in \PP_2(\CC)$ be Hermitian cubic over $K$, and
let $z'_0, z''_0$ be its other two Galois conjugates. Let
$\underline{G}$ be the linear algebraic group defined over $\QQ$, such
that $\underline{G}(\ZZ)=\PSU_q(\OOO_K)$ and $\underline{G}(\RR) =
\PSU_q$. It has $\RR$-rank one. The pointwise stabiliser of $\{z'_0,
z''_0\}$ in $\underline{G}$ is the centraliser $Z(\underline{T})$ of a
maximal algebraic torus $\underline{T}$ in $\underline{G}$, since
$z'_0, z''_0$ are distinct and isotropic. Note that $Z(\underline{T})$
also fixes $z_0$, since $z_0$, being orthogonal to $z'_0$ and $z''_0$,
belongs to the complex projective lines tangent to the null cone of
$q$ in $\PP_2(\CC)$ at $z'_0$ and $z''_0$, and as seen above, these
two projective lines meet at exactly one point. The algebraic group
$Z(\underline{T})$, being the pointwise stabiliser of $\{z_0, z'_0,
z''_0\}$ which is invariant under $\Gal(\CC/K)$, is defined over $K$.
The torus $\underline{T}$ has rank one and has finite index in
$Z(\underline{T})$. Hence $\underline{T}$ is also defined over $K$,
and is isomorphic to $\CC^\times$ over $\CC$. It has no nontrivial
$\QQ$-character (as it is one-dimensional, it is not defined over
$\QQ$, otherwise its fixed points $z_0,z'_0, z''_0$ would individually
be defined over $\QQ$). Hence by the Borel and Harish-Chandra theorem
(see \cite[Th.~9.4]{BorHar62}, though the particular case we use here
is due to Ono), $\underline{T}(\ZZ)$ is a lattice in $\underline{T}
(\RR)$.  Such a lattice contains an element of infinite order
$\ga_0$. The set of fixed points of $\ga_0$ in $\PP_2(\CC)$ is
$\{z_0,z'_0,z''_0\}$. Any proper nonzero linear subspace of $\CC^3$
invariant under $\ga_0$ is the sum of one or two geodesic lines in
$\{z_0,z'_0,z''_0\}$, hence is not defined over $K$ since $z_0$ is
cubic. Therefore $\ga_0$ is $K$-irreducible. As seen above, $z_0$ is
the only fixed point of $\ga_0$ in the positive cone of $q$. \cqfd

\medskip Let $\ga_0$ be as in Proposition \ref{prop:caracfixpoint} for
$z_0$ given by the statement of Theorem \ref{theo:countingcubic}. As
seen in the above proof, $\ga_0$ is loxodromic.  Let $D^+$ be the
geodesic line in the projective model of $\hdc$ with endpoints the
other two Galois conjugates $z'_0,z''_0$ of $z_0$ over $K$, and let
$G_{D^+}$ be its stabiliser in $G$. Up to replacing $\ga_0$ by a
power, we may assume that $\ga_0\in G$. We may also assume that
$\ga_0$ is primitive in $G$, so that if $\lambda_0$ is its eigenvalue
with modulus $>1$, then its translation length in $\hdc$ is
$\ln|\lambda_0|$ (see Equation \eqref{eq:calctransllengthloxo}), and
\begin{equation}\label{eq:covolGDplus}
\Vol(G_{D^+}\bs D^+)=\frac{\ln|\lambda_0|}{\iota_0}\,,
\end{equation}
with $\iota_0$ defined in the statement of Theorem
\ref{theo:countingcubic}.

Let $G_{z_0}$ be the stabiliser of $z_0$ in $G$, which is also the
stabiliser of $z_0^\perp=D^+$, hence coincides with $G_{D^+}$. Let
$g\in G$ be such that the geodesic line $gD^+$ is disjoint from $\H_1$
(which is the case except for finitely many double classes in
$G_{\H_1}\bs G/G_{D^+}$). Let $z'$, $z''$ be the endpoints of $gD^+$.
Let $\delta_g$ be the common perpendicular from $\H_1$ to $gD^+$. The
length of $\delta_g$ is, by \cite[Lem.~3.4]{ParPau11MZ}, by Equation
\eqref{eq:distentrhorob}, and by the definition of the complexity,
$$
\ell(\delta_g)=d(gD^+,\H_2)-d(\H_1,\H_2)=
-\ln\frac{d''_{\rm Cyg}(z',z'')}{2}-\frac{\ln 2}{2}=\ln(\sqrt{2}\,c(gz_0))\;.
$$
Therefore, by Corollary \ref{coro:complexhyperbo} (with $n=2$, in the
case $D^-=\H_1$ is a horoball, whose pointwise stabiliser is trivial,
and $D^+$ is a geodesic line, with pointwise stabiliser of order $n_0$
as defined in the statement of Theorem \ref{theo:countingcubic}), and
by Equations \eqref{eq:covolG}, \eqref{eq:covolGinfty} and
\eqref{eq:covolGDplus},
\begin{align*}
&\card\{z\in G_\infty\bs G\cdot z_0\;:\; c(z)\leq s\} =
\card\{[g]\in G_\infty\bs G/G_{z_0}\;:\; c(gz_0)\leq s\}\\ =\;&
\card\{[g]\in G_{D^-}\bs G/G_{D^+}\;:\; \ell(\delta_g)\leq 
\ln (s\,\sqrt{2}\,)\}+\bigO(1)\\ =\;&
\N_{D^-,\,D^+}\big(\ln (s\,\sqrt{2}\,)\big)+ \bigO(1)=
\frac{c(D^-,D^+)}{n_0}\;e^{4\ln (s\,\sqrt{2}\,)}\;
(1+\bigO(e^{-\kappa\ln (s\,\sqrt{2}\,)}))\\ =\;&
\frac{32\,\Vol(G_{D^-}\bs D^-)\Vol(G_{D^+}\bs D^+)}
{9\,\pi\,n_0\,\Vol(G\bs \hdc)}\;s^4\,(1+\bigO(s^{-\kappa}))\\ =\;&
\frac{256\,\ln|\lambda_0|\;\zeta(3)\;n_ \infty}
{3\,\iota_0\,|\OOO_K^\times|\,|D_K|^{3/2}\,\zeta_K(3)\,n_0\,
[\,\overline{\Ga}\,:G]}\;s^4\,(1+\bigO(s^{-\kappa}))\,.
\end{align*}
This proves Theorem \ref{theo:countingcubic}. \cqfd


\section{Counting arithmetic chains in hyperspherical 
geometry}
\label{sect:chains}

Let us consider again the Hermitian form $q=-z_0\,\overline{z_2}
-z_2\,\overline{z_0} + |z_1|^2$ of signature $(1,2)$ on $\CC^3$ with
coordinates $(z_0,z_1,z_2)$. Following Poincaré \cite{Poincare07} (who
was rather using the diagonal form $-|z_0|^2 + |z_1|^2 +|z_2|^2$), see
also \cite{Cartan32}, we will call {\it hypersphere} the projective
isotropic locus of $q$, that is the subspace
$$
\HS=\{[z_0:z_1:z_2]\in \PP_2(\CC)\;:\; q(z_0,z_1,z_2)=0\}
$$
of the complex projective plane $\PP_2(\CC)$ with homogeneous
coordinates $[z_0:z_1:z_2]$. It is a real analytic submanifold,
diffeomorphic to the $3$-sphere $\SSS_3$. The subgroup $\PSU_q$ of
$\PSL_3(\CC)$, acting projectively on $\PP_2(\CC)$, preserves the
hypersphere.

In Section \ref{sec:cxhyp}, we introduced a natural modification
$d''_{\rm Cyg}$ of Cygan's distance $d_{\rm Cyg}$ on $\HS-\{\infty\}$,
with $\infty =[1:0:0]$, as follows. We identified $\HS-\{\infty\}$
with the real quadric (called a hyperconic by Segre) $\Heis_3=
\{(w_0,w) \in \CC^2 \;:\; 2\,\Re\; w_0 -|w|^2=0\}$ by the map $(w_0,w)
\mapsto [w_0:w:1]$. Then $d''_{\rm Cyg}$ is the unique map from
$(\HS-\{\infty\})^2$ to $[0,+\infty[$, invariant under the diagonal
action of the unipotent radical of the stabiliser of $\infty$ in
$\PSU_q$, such that
$$ 
d''_{\rm Cyg}((w_0,w),(0,0))^2=
\frac{|w|^4+4\,\Im^2\;w_0}{(|w|^4+4\,\Im^2\;w_0)^{\frac{1}{2}}+|w|^2}=
\frac{4\,|w_0|^2}{2\,|w_0|+|w|^2}\;.
$$

A complex projective line $L$ in $\PP_2(\CC)$ intersects the
hypersphere either in the empty set, or a one point set (in which case
$L$ is the unique complex projective line tangent to $\HS$ at this
point, giving at this point the canonical contact structure of $\HS$),
or a real analytic circle, called a {\it chain} (a notion attributed
to von Staudt by \cite[footnote 3)]{Cartan32}).  A chain $C$ separates
the complex projective line $L(C)$ containing it into two real discs
$D_\pm(C)$, which we endow with their unique Poincaré metric (of
constant curvature $-1$) invariant under the stabiliser of $C$ in
$\PSU_q$.  It is finite if and only if $C$ is a finite chain, that is,
if it does not contain $\infty=[1:0:0]$.  We refer to \cite[\S
4.3]{Goldman99} for more informations on the chains, including the
following fact: the infinite chains are precisely the fibers of the
{\it vertical projection} $(w_0,w)\mapsto w$, and the finite chains
are ellipses in the affine coordinates $(2\,\Im\;w_0,w)$ of $\Heis_3$
whose vertical projections are circles.

\begin{center}
\includegraphics[scale=.8]{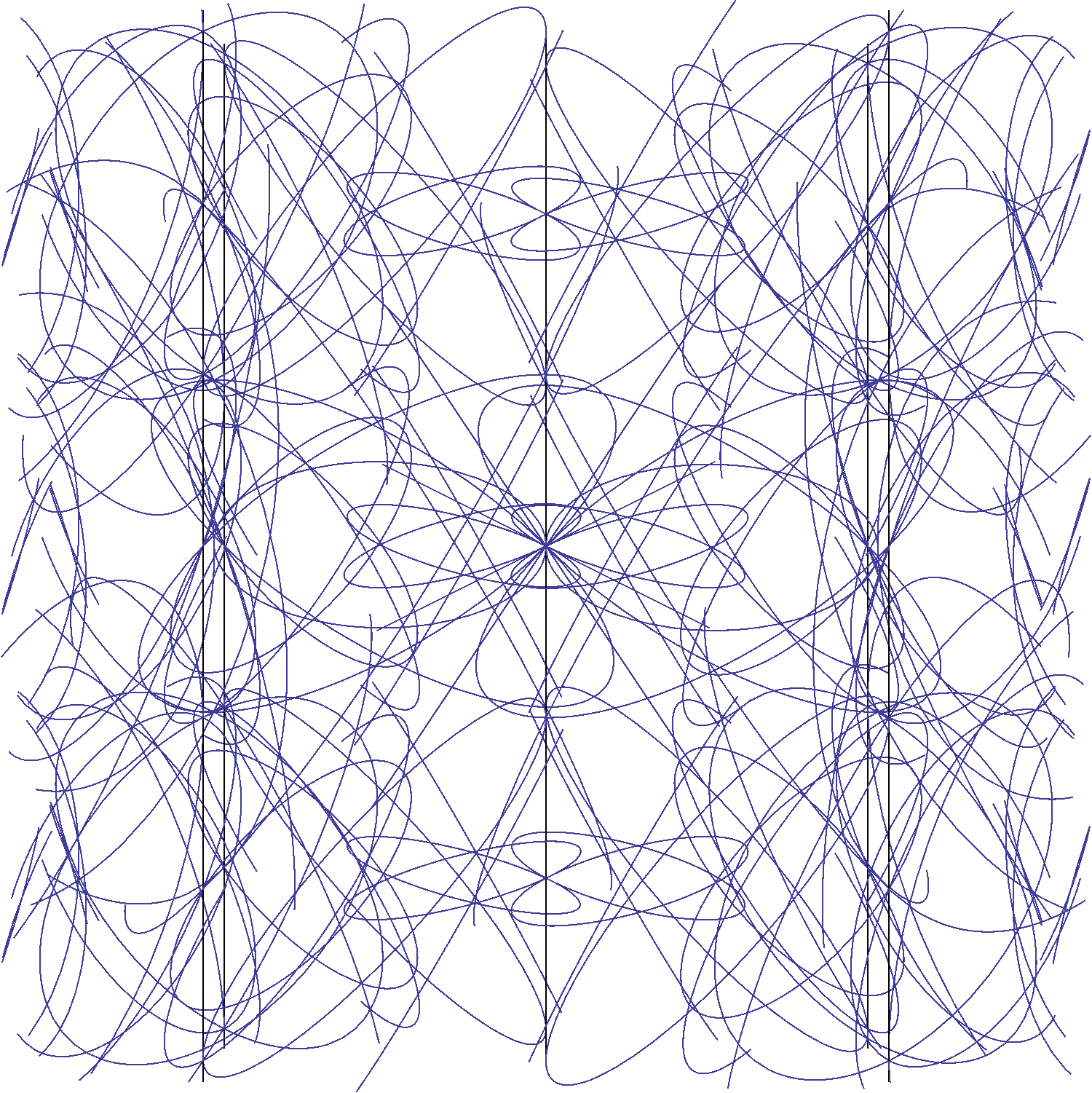}
\end{center}

Let $K$ be an imaginary quadratic number field. For every finite index
subgroup $G$ of the arithmetic lattice $\PSU_q(\OOO_K)$, we will
denote by $G_C$ the stabiliser of $C$ in $G$, by $G_\infty$ the
stabiliser of $\infty$ in $G$, and by $\covol_G(C)$ the (common)
volume of $G_C\bs D_\pm(C)$.  A chain $C$ will be called {\it
  arithmetic} (over $K$) if $\PSU_q(\OOO_K)_C$ has a dense orbit in
$C$.

\btheo \label{theo:countchain} Let $C_0$ be an arithmetic chain in the
hypersphere $\HS$ over an imaginary quadratic number field
$K$. Let $G$ be a finite index subgroup of $\PSU_q(\OOO_K)$. Then
there exists a constant $\kappa>0$ such that, as $\epsilon>0$ tends to
$0$, the number $\psi_{C_0,\,G}(\epsilon)$ of chains modulo $G_\infty$
in the $G$-orbit of $C_0$ with $d_{\rm Cyg}$-diameter at least
$\epsilon$ is equal to
$$
\frac{2048\,\zeta(3)\,\covol_G(C_0)\,[\PSU_q(\OOO_K)_\infty:G_\infty]}
{|\OOO_K^\times|\,|D_K|^{\frac{3}{2}}\,\zeta_K(3)\,n_{0,\,G}\,[\PSU_q(\OOO_K):G]}
\;\epsilon^{-4}\big(1+\bigO(\epsilon^\kappa)\big)\;,
$$
where $n_{0,\,G}$ is the order of the pointwise stabiliser of $C_0$ in $G$.
\etheo

\medskip Given a complex projective line $L$ in $\PP_2(\CC)$, there is
a unique order $2$ complex projective map with fixed point set $L$,
called the {\it reflexion} in $L$. Given a finite chain $C$, contained
in the projective line $L(C)$, the {\it center} of $C$ (see for
instance \cite[4.3.3]{Goldman99}), denoted by $\operatorname{cen}(C)\in
\HS-\{\infty\}=\Heis_3$, is the image of $\infty=[1:0:0]$ under
the reflexion in $L(C)$. The following result is an equidistribution
result in the Heisenberg group of the centers of the arithmetic chains
in a given orbit under (a finite index subgroup of) $\PSU_q(\OOO_K)$.

\btheo\label{theo:equidischain} Let $C_0$ and $G$ be as in Theorem
\ref{theo:countchain}. As $\epsilon>0$ tends to $0$, we have
$$
\frac{n_{0,\,G}\,(1+2\,\delta_{D_K,-3})
  \,|D_K|^{\frac{5}{2}}\,\zeta_K(3)\,[\PSU_q(\OOO_K):G]}
{1024\,\zeta(3)\,\covol_G(C_0)}\;\epsilon^{4}
\sum_{C\in G\cdot C_0\;:\;\diam_{d_{\rm Cyg}}\;C\geq \epsilon}\;
\Delta_{\operatorname{cen}(C)}\;\weakstar\;\haarheis\,.
$$ 
\etheo

\noindent{\bf Proof of Theorem \ref{theo:countchain} and 
Theorem \ref{theo:equidischain}. } 
As seen in Section \ref{sec:cxhyp}, the hypersphere $\HS$ is the
boundary at infinity of the projective model of the complex hyperbolic
space $\hdc$. The chains are precisely the boundary at infinity of the
complex geodesic lines in $\hdc$. The diameter of a chain is invariant
under the stabiliser in $\PSU_q$ of the horosphere $\partial \H_1$,
hence is invariant under $G_\infty$. The counting function $\psi_{C_0,\,G}$
is thus well defined.

Recall (see for instance \cite{Bowditch95}, in particular for the
terminology) that a discrete group of isometries $\Ga$ of a complete
simply connected Riemannian manifold $M$ with sectional curvature at
most $-1$ is {\it geometrically finite} if and only if every limit
point of $\Ga$ is either a bounded parabolic point or a conical limit
point.  Furthermore, the discrete groups $\Ga$ of isometries of $M$
with finite covolume are the geometrically finite discrete groups of
isometries $\Ga$ whose limit set is the whole sphere at infinity
$\partial_\infty M$ of $M$; then the orbit under $\Ga$ of every point
in $\partial_\infty M$ is dense in $\partial_\infty M$.

Let $C$ be a chain in $\HS$, and let $D$ be the complex geodesic
line (which is totally geodesic in $\hdc$) with $\partial_\infty
D=C$. Hence, $C$ is arithmetic over $K$ if and only if the stabiliser
in $\PSL_q(\OOO_K)$ (or equivalently in $G$) of $D$ has finite
covolume on $D$.

We denote by $D^+$ the complex geodesic line in $\hdc$ with
$\partial_\infty D^+=C_0$. Let $G_{D^+}$ be the stabiliser of $D^+$ in
$G$. By definition, we have
\begin{equation}\label{eq:covolVol}
\Vol(G_{D^+}\bs D^+)=4\;\covol_G(C_0)\,,
\end{equation}
since the sectional curvature of $D^+$ is constant $-4$ and $D^+$
has real dimension $2$.

\medskip Let $g\in G$ be such that the complex geodesic line $gD^+$ is
disjoint from $\H_1$ (which is the case except for $g$ in finitely
many double classes in $G_{\H_1}\bs G/G_{D^+}$). Let $\delta_g$ be the
common perpendicular from $\H_1$ to $gD^+$. Its length
$\ell(\delta_g)$ is the minimum of the distances from $\H_1$ to a
geodesic line between two points of $\partial_\infty(gD^+)=g
C_0$. Hence, by \cite[Lem.~3.4]{ParPau11MZ}, Equation
\eqref{eq:distentrhorob} and Lemma \ref{lem:diamchain}, we have
\begin{align}\label{eq:relatlongcomperpchain}
\ell(\delta_g)= &
\min_{x,y\in gC_0,\;x\neq y}\ln \frac{2}{d''_{\rm Cyg}(x,y)} - d(\H_1,\H_2)
=-\max_{x,y\in gC_0,\;x\neq y}\ln \frac{d''_{\rm Cyg}(x,y)}{2}-
\frac{\ln 2}{2}\nonumber\\
= & 
-\ln \frac{\diam_{d''_{\rm Cyg}}(gC_0)}{\sqrt{2}}=
-\ln \frac{\diam_{d_{\rm Cyg}}(gC_0)}{2}\,.
\end{align}

Now, we apply Corollary \ref{coro:complexhyperbo} with $n=2$, in the
case $D^-=\H_1$ is a horoball, whose pointwise stabiliser is trivial,
and $D^+$ is a complex geodesic line, whose pointwise stabiliser has
order $n_{0,\,G}$ as defined in the statement of Theorem
\ref{theo:countchain}. Respectively by the definition of the counting
function $\psi_{C_0,\,G}$ in the statement of Theorem
\ref{theo:countchain}, since the stabiliser of $C_0$ in $G$ is equal
to $G_{D^+}$, by Equation \eqref{eq:relatlongcomperpchain}, by
Corollary \ref{coro:complexhyperbo}, and by Equations
\eqref{eq:covolG}, \eqref{eq:covolGinfty} and \eqref{eq:covolVol}, we
have, as $\epsilon>0$ tends to $0$,
\begin{align*}
\psi_{C_0,\,G}(\epsilon)& 
=\card\{C\in G_\infty\bs G\cdot C_0\;:\; 
\diam_{d_{\rm Cyg}}(C)\geq \epsilon\}\\ & =
\card\{[g]\in G_\infty\bs G/G_{D^+}\;:\; 
\diam_{d_{\rm Cyg}}(gC_0)\geq \epsilon\}\\ & =
\card\{[g]\in G_{\H_1}\bs G/G_{D^+}\;:\; \ell(\delta_g)
\leq -\ln\frac{\epsilon}{2}\} +\bigO(1)\\ & =
\N_{\H_1,\,D^+}(-\ln\frac{\epsilon}{2})+\bigO(1)=
\frac{c(\H_1,D^+)}{n_{0,\,G}}\;e^{-4\ln\frac{\epsilon}{2}}
\big(1+\bigO(e^{\kappa \ln\frac{\epsilon}{2}})\big)\\ & =
\frac{64\,\Vol(G_{\H_1}\bs \H_1)\,\Vol(G_{D^+}\bs D^+)}
{3\,\pi\,n_{0,\,G}\,\Vol(G\bs \hdc)}\;\epsilon^{-4}
\big(1+\bigO(\epsilon^\kappa)\big)\\ & =
\frac{2048\,\zeta(3)\,\covol_G(C_0)\,[\PSU_q(\OOO_K)_\infty:G_\infty]}
{|\OOO_K^\times|\,|D_K|^{\frac{3}{2}}\,\zeta_K(3)\,n_{0,\,G}\,[\PSU_q(\OOO_K):G]}
\;\epsilon^{-4}\big(1+\bigO(\epsilon^\kappa)\big)\;.
\end{align*}

This proves Theorem \ref{theo:countchain}. Let us now prove Theorem
\ref{theo:equidischain}.

\medskip We apply the equidistribution result in Equation
\eqref{eq:distribhorobcomplexhyp} (in the case $D^+$ is a complex
geodesic line) of the origins $\operatorname{or}(\delta_g)$ of the common
perpendiculars $\delta_g$ from $D^-=\H_1$ to the images $gD^+$ for
$g\in G$. As $t\ra+\infty$, we hence have
\begin{equation}\label{eq:distribcentchain}
\frac{3\, n_{0,\,G}\,\pi\,\Vol(G\bs \hdc)}{\Vol(G_{D^+}\bs D^+)}
\;e^{-4\,t}\;\sum_{[g]\in G/G_{D^+},\; \ell(\delta_g)\leq t}\;
\Delta_{\operatorname{or}(\delta_g)}\;\weakstar\; \Vol_{\partial \H_1}\,.
\end{equation}

Let $f:\partial_\infty\hdc-\{\infty\}=\Heis_3\ra \partial \H_1$ be the
orthogonal projection map, which in horospherical coordinates is
$(\zeta,u,0)\mapsto (\zeta,u,1)$. By Equations \eqref{eq:volHun} and
\eqref{eq:relatlambdatroihaarheis}, the image of the Haar measure
$\haarheis$ by $f$ is
$$
f_*\haarheis=\vol_{\partial \H_1}\,.
$$ 

Note that, for every chain $C$, if $r_C$ is the reflexion on the
complex projective line containing $C$, then the geodesic line from
$\infty$ to $\operatorname{cen}(C)=r_C(\infty)$, being invariant under
$r_C$, is orthogonal to the complex geodesic line with boundary at
infinity $C$.  Hence for every $g\in G$, we have
$$
f^{-1}(\operatorname{or}(\delta_g))=\operatorname{cen}(gC_0)\;.
$$

Let us use in Equation \eqref{eq:distribcentchain} the change of
variables $t=-\ln\,\frac{\epsilon}{2}$ and the continuity of the
pushforward of measures by $f^{-1}$. By Equations \eqref{eq:covolG},
\eqref{eq:covolVol} and \eqref{eq:relatlongcomperpchain}, as
$\epsilon>0$ tends to $0$, the measures
$$
 \frac{n_{0,\,G}\,(1+2\,\delta_{D_K,-3})\,|D_K|^{\frac{5}{2}}\,\zeta_K(3)\,
[\PSU_q(\OOO_K):G]}{1024\,\zeta(3)\,\covol_G(C_0)}\;\epsilon^{4}
\sum_{[g]\in  G/G_{D^+},\; \diam_{d_{\rm Cyg}}(gC_0)\geq \epsilon}\; 
\Delta_{\operatorname{cen}(gC_0)}
$$
weak-star converge to the Haar measure $\haarheis$. 
This proves Theorem \ref{theo:equidischain}. \cqfd

\medskip
\noindent{\bf Example.} Let $K=\QQ(i)$, so that $|D_K|=|\OOO_K^\times|
=4$. Let $C_0$ be the intersection of the hypersphere and the complex
projective line $L(C_0)$ with equation $z_1=0$.  

Let us use $\infty= [1:0:0]$ as the point at infinity in $L(C_0)$, and
let us identify $L(C_0)$ with $\PP_1(\CC)=\CC\cup\{\infty\}$ by the
map $[z_0:0:z_2]\mapsto z=\frac{z_0}{z_2}$. Then the restriction
$\Ga_{1,\,1}$ to $L(C_0)$ of the stabiliser in $\PSU_q(\OOO_K)$ of
$C_0$ is exactly the subgroup of $\PSL_2(\CC)$ preserving the right
halfplane with equation $\Re \;z>0$. The homography $z\mapsto iz$
sends the right halfplane to the upper halfplane $\HH^2_\RR$. Hence
$\begin{pmatrix} i & 0\\ 0& 1\end{pmatrix}\Ga_{1,\,1}\begin{pmatrix}
  -i & 0\\ 0& 1\end{pmatrix}$ is the subgroup of $\PSL_2(\CC)$
preserving $\HH^2_\RR$ and having coefficients in $\ZZ[i]$. It is
hence equal to $\PSL_2(\ZZ)$, and it is well known that the hyperbolic
volume of $\PSL_2(\ZZ)\bs\HH^2_\RR$ is $\frac{\pi}{3}$. In particular,
$C_0$ is an arithmetic chain over $K$.

Let $\alpha$ be an element in the pointwise stabiliser of $C_0$ in
$\PSU_q(\OOO_K)$. Then $\alpha$ fixes $\infty=[1:0:0]$ and $[0:0:1]$,
hence consists of diagonal matrices. The diagonal coefficients
$\lambda_1, \lambda_2,\lambda_3$ of $\alpha$ belong to $\ZZ[i]^\times=
\{\pm 1,\pm i\}$. We have $\lambda_1=\lambda_3$, otherwise $\alpha$
would not act by the identity on $C_0$. Therefore $(\lambda_1,
\lambda_2,\lambda_3)$ belongs to $\{(1,1,1),(-1,1,-1),(i,-1,i),$
$(-i,-1,-i)\}$.  This gives the following values
$$
\covol_{\PSU_q(\ZZ[i])}(C_0) = \frac{\pi}{3}\;\;\;{\rm and}
\;\;\;n_{0,\,\PSU_q(\ZZ[i])}= 4\;.
$$
Theorems \ref{theo:countchain} and \ref{theo:equidischain} then give
$$
\psi_{C_0,\,\PSU_q(\ZZ[i])}(\epsilon)=
\frac{16\,\pi\,\zeta(3)}{3\,\zeta_{\QQ(i)}(3)} 
\;\epsilon^{-4}\big(1+\bigO(\epsilon^\kappa)\big)\;, 
$$
and
$$
\frac{3\,\zeta_{\QQ(i)}(3)}{4\,\pi\,\zeta(3)}\;\epsilon^{4}
\sum_{C\in \PSU_q(\ZZ[i])\cdot C_0\;:\;\diam_{d_{\rm Cyg}}\;C\geq \epsilon}\;
\Delta_{\operatorname{cen}(C)}\;\weakstar\;\haarheis\,.
$$ 

The above picture, as well as the one in the introduction, represents
the orbit of $C_0$ under the arithmetic lattice $\PSU_q(\ZZ[i])$. The
pictures were produced using Mathematica and they show images of the
chain $C_0$ by elements of $\PSU_q(\ZZ[i])$ of word length at most 16
in the generating set given in \cite[Theo.~7]{FalFraPar11}. Both
images show the same set of chains from two different viewpoints in
the $3$-dimensional space $\Heis_3$, except that in the second one we
removed some of the smaller chains in order to make the structure more visible.

{\small \bibliography{../biblio} }

\bigskip
{\small
\noindent \begin{tabular}{l} 
Department of Mathematics and Statistics, P.O. Box 35\\ 
40014 University of Jyv\"askyl\"a, FINLAND.\\
{\it e-mail: jouni.t.parkkonen@jyu.fi}
\end{tabular}
\medskip

\noindent \begin{tabular}{l}
D\'epartement de math\'ematique, UMR 8628 CNRS, B\^at.~425\\
Universit\'e Paris-Sud,
91405 ORSAY Cedex, FRANCE\\
{\it e-mail: frederic.paulin@math.u-psud.fr}
\end{tabular}
}

\end{document}